# Computational techniques for proving identities in symmetric compositions

*Pablo Alberca Bjerregaard[a] and Cándido Martín González[b]*

[a] University of Málaga, Department of Applied Mathematics. Málaga. Spain. E-mail: pgalberca@uma.es

[b] University of Málaga, Department of Algebra, Geometry and Topology. Málaga. Spain. E-mail: candido@apncs.cie.uma.es

January - 2007

**Abstract**. We present in this work a complete session in a *Mathematica* notebook. The aim of this notebook is to check identities in symmetric compositions. This notebook is a complement of our work [1] and it has all the explicit computations. We refer the reader to that paper which can be seen in `http://www.uibk.ac.at/mathematik/loos/jordan/index.html`. First of all we will present a few number of comands in order to simplify identities by extracting scalars, **SOut**. The rest of the strategy holds on the powerfull of using patterns and rules.

## ■ Rules and functions for linearization

```
expandq =
  {q[x_ + y_] → b[x, y] + q[x] + q[y], q[x_ · y_] → q[x] q[y]};

expandb = {b[x_ + y_, z_] → b[x, z] + b[y, z],
    b[x_, y_ + z_] → b[x, y] + b[x, z],
    b[-x_, y_] → -b[x, y], b[x_, - y_ · z_] → -b[x, y · z],
    b[x_, -(z_ · t_)] → -b[x, z · t]};

expanddot = {(x_ + y_) · z_ → x · z + y · z,
    x_ · (y_ + z_) → x · y + x · z, x_ · -y_ → - x · y};

SOut[b[x_, y_], scal_List] :=
  Module[{scaln = Union[scal, Table[ξ_i, {i, Length[scal]}]],
    m = Length[scal], res = b[x, y]},
   Do[res = res //. -scal[[j]] → ξ_j, {j, m}];
   Do[
    If[Not[FreeQ[res, scaln[[j]]]],
     res = scaln[[j]]^Length[Position[res,scaln[[j]]]]
       (res /. scaln[[j]] → 1)], {j, Length[scaln]}];
   Do[res = res //. ξ_j → -scal[[j]], {j, m}];
   res];
```

```
SOut[x_CenterDot, scal_List] :=
  Module[{scaln = Union[scal, Table[ξᵢ, {i, Length[scal]}]],
    m = Length[scal], res = x},
   Do[res = res //. -scal[[j]] → ξⱼ, {j, m}];
   Do[
    If[Not[FreeQ[res, scaln[[j]]]],
     res = scaln[[j]]^Length[Position[res,scaln[[j]]]]
        (res /. scaln[[j]] → 1)], {j, Length[scaln]}];
   Do[res = res //. ξⱼ → -scal[[j]], {j, m}];
   res];

SOut[q[α_. * x_], scal_List] := Which[MemberQ[scal, α],
  α^2 q[x], MemberQ[scal, x], x^2 q[α], True, q[α x]]

SOut[x_Symbol, scal_List] := x; SOut[x_Power, scal_List] := x;
SOut[x_Subscript, scal_List] := x;
SOut[x_?NumberQ, scal_List] := x;

SOut[x_Times, scal_List] := Map[SOut[#, scal] &, x];
SOut[x_Plus, scal_List] := Map[SOut[#, scal] &, x]
```

## ■ Linearization of $q(x \cdot y) = q(x) \, q(y)$

```
comp = q[x·y] - q[x] q[y];

comp //. {x → z₁ + α z₂, y → z₃ + β z₄}
```

$q((z_1 + \alpha z_2) \cdot (z_3 + \beta z_4)) - q(z_1 + \alpha z_2) \, q(z_3 + \beta z_4)$

```
% //. expanddot
```

$q(z_1 \cdot z_3 + z_1 \cdot (\beta z_4) + (\alpha z_2) \cdot z_3 + (\alpha z_2) \cdot (\beta z_4)) - q(z_1 + \alpha z_2) \, q(z_3 + \beta z_4)$

```
% //. expandq
```

$b(z_1 \cdot z_3, z_1 \cdot (\beta z_4) + (\alpha z_2) \cdot z_3 + (\alpha z_2) \cdot (\beta z_4)) + b(z_1 \cdot (\beta z_4), (\alpha z_2) \cdot z_3 + (\alpha z_2) \cdot (\beta z_4)) +$
$b((\alpha z_2) \cdot z_3, (\alpha z_2) \cdot (\beta z_4)) + q(z_1) \, q(z_3) + q(\alpha z_2) \, q(z_3) + q(z_1) \, q(\beta z_4) +$
$q(\alpha z_2) \, q(\beta z_4) - (b(z_1, \alpha z_2) + q(z_1) + q(\alpha z_2))(b(z_3, \beta z_4) + q(z_3) + q(\beta z_4))$

```
% //. expandb
```

$b(z_1 \cdot z_3, z_1 \cdot (\beta z_4)) + b(z_1 \cdot z_3, (\alpha z_2) \cdot z_3) + b(z_1 \cdot z_3, (\alpha z_2) \cdot (\beta z_4)) +$
$b(z_1 \cdot (\beta z_4), (\alpha z_2) \cdot z_3) + b(z_1 \cdot (\beta z_4), (\alpha z_2) \cdot (\beta z_4)) +$
$b((\alpha z_2) \cdot z_3, (\alpha z_2) \cdot (\beta z_4)) + q(z_1) \, q(z_3) + q(\alpha z_2) \, q(z_3) + q(z_1) \, q(\beta z_4) +$
$q(\alpha z_2) \, q(\beta z_4) - (b(z_1, \alpha z_2) + q(z_1) + q(\alpha z_2))(b(z_3, \beta z_4) + q(z_3) + q(\beta z_4))$

```
Expand[%]
```

$b(z_1 \cdot z_3, z_1 \cdot (\beta z_4)) + b(z_1 \cdot z_3, (\alpha z_2) \cdot z_3) + b(z_1 \cdot z_3, (\alpha z_2) \cdot (\beta z_4)) +$
$b(z_1 \cdot (\beta z_4), (\alpha z_2) \cdot z_3) + b(z_1 \cdot (\beta z_4), (\alpha z_2) \cdot (\beta z_4)) +$
$b((\alpha z_2) \cdot z_3, (\alpha z_2) \cdot (\beta z_4)) - b(z_1, \alpha z_2) \, b(z_3, \beta z_4) - b(z_3, \beta z_4) \, q(z_1) -$
$b(z_3, \beta z_4) \, q(\alpha z_2) - b(z_1, \alpha z_2) \, q(z_3) - b(z_1, \alpha z_2) \, q(\beta z_4)$

```
SOut[%, {α, β}]
```

$\beta\, b(z_2 \cdot z_3, z_2 \cdot z_4)\, \alpha^2 - \beta\, b(z_3, z_4)\, q(z_2)\, \alpha^2 + b(z_1 \cdot z_3, z_2 \cdot z_3)\, \alpha + \beta\, b(z_1 \cdot z_3, z_2 \cdot z_4)\, \alpha +$
$\beta\, b(z_1 \cdot z_4, z_2 \cdot z_3)\, \alpha + \beta^2\, b(z_1 \cdot z_4, z_2 \cdot z_4)\, \alpha - \beta\, b(z_1, z_2)\, b(z_3, z_4)\, \alpha -$
$b(z_1, z_2)\, q(z_3)\, \alpha - \beta^2\, b(z_1, z_2)\, q(z_4)\, \alpha + b(z_1 \cdot z_3, z_1 \cdot z_4) - \beta\, b(z_3, z_4)\, q(z_1)$

```
CoefficientList[%, {α, β}]
```

$$\begin{pmatrix} 0 & b(z_1 \cdot z_3, z_1 \cdot z_4) - b(z_3, z_4)\, q(z_1) & 0 \\ b(z_1 \cdot z_3, z_2 \cdot z_3) - b(z_1, z_2)\, q(z_3) & b(z_1 \cdot z_3, z_2 \cdot z_4) + b(z_1 \cdot z_4, z_2 \cdot z_3) - b(z_1, z_2)\, b(z_3, z_4) & b(z_1 \cdot z_4, z_2 \cdot z_4) \\ 0 & b(z_2 \cdot z_3, z_2 \cdot z_4) - b(z_3, z_4)\, q(z_2) & 0 \end{pmatrix}$$

## ■ Linearization of $(x \cdot y) \cdot x = q(x)\, y$

```
iden = (x · y) · x - q[x] y
```

$(x \cdot y) \cdot x - y\, q(x)$

```
iden //. {x → z₁ + α z₂, y → z₃ + β z₄}
```

$((z_1 + \alpha z_2) \cdot (z_3 + \beta z_4)) \cdot (z_1 + \alpha z_2) - q(z_1 + \alpha z_2)\, (z_3 + \beta z_4)$

```
% //. expanddot
```

$(z_1 \cdot z_3) \cdot z_1 + (z_1 \cdot z_3) \cdot (\alpha z_2) + (z_1 \cdot (\beta z_4)) \cdot z_1 + (z_1 \cdot (\beta z_4)) \cdot (\alpha z_2) + ((\alpha z_2) \cdot z_3) \cdot z_1 +$
$((\alpha z_2) \cdot z_3) \cdot (\alpha z_2) + ((\alpha z_2) \cdot (\beta z_4)) \cdot z_1 + ((\alpha z_2) \cdot (\beta z_4)) \cdot (\alpha z_2) - q(z_1 + \alpha z_2)\, (z_3 + \beta z_4)$

```
% //. expandq
```

$(z_1 \cdot z_3) \cdot z_1 + (z_1 \cdot z_3) \cdot (\alpha z_2) + (z_1 \cdot (\beta z_4)) \cdot z_1 + (z_1 \cdot (\beta z_4)) \cdot (\alpha z_2) +$
$((\alpha z_2) \cdot z_3) \cdot z_1 + ((\alpha z_2) \cdot z_3) \cdot (\alpha z_2) + ((\alpha z_2) \cdot (\beta z_4)) \cdot z_1 +$
$((\alpha z_2) \cdot (\beta z_4)) \cdot (\alpha z_2) - (b(z_1, \alpha z_2) + q(z_1) + q(\alpha z_2))\, (z_3 + \beta z_4)$

```
Expand[%]
```

$(z_1 \cdot z_3) \cdot z_1 + (z_1 \cdot z_3) \cdot (\alpha z_2) + (z_1 \cdot (\beta z_4)) \cdot z_1 + (z_1 \cdot (\beta z_4)) \cdot (\alpha z_2) +$
$((\alpha z_2) \cdot z_3) \cdot z_1 + ((\alpha z_2) \cdot z_3) \cdot (\alpha z_2) + ((\alpha z_2) \cdot (\beta z_4)) \cdot z_1 + ((\alpha z_2) \cdot (\beta z_4)) \cdot (\alpha z_2) -$
$b(z_1, \alpha z_2)\, z_3 - q(z_1)\, z_3 - q(\alpha z_2)\, z_3 - \beta\, b(z_1, \alpha z_2)\, z_4 - \beta\, q(z_1)\, z_4 - \beta\, q(\alpha z_2)\, z_4$

```
SOut[%, {α, β}]
```

$(z_2 \cdot z_3) \cdot z_2\, \alpha^2 + \beta\, (z_2 \cdot z_4) \cdot z_2\, \alpha^2 - q(z_2)\, z_3\, \alpha^2 - \beta\, q(z_2)\, z_4\, \alpha^2 + (z_1 \cdot z_3) \cdot z_2\, \alpha +$
$\beta\, (z_1 \cdot z_4) \cdot z_2\, \alpha + (z_2 \cdot z_3) \cdot z_1\, \alpha + \beta\, (z_2 \cdot z_4) \cdot z_1\, \alpha - b(z_1, z_2)\, z_3\, \alpha -$
$\beta\, b(z_1, z_2)\, z_4\, \alpha + (z_1 \cdot z_3) \cdot z_1 + \beta\, (z_1 \cdot z_4) \cdot z_1 - q(z_1)\, z_3 - \beta\, q(z_1)\, z_4$

```
CoefficientList[%, {α, β}]
```

$$\begin{pmatrix} (z_1 \cdot z_3) \cdot z_1 - q(z_1)\, z_3 & (z_1 \cdot z_4) \cdot z_1 - q(z_1)\, z_4 \\ (z_1 \cdot z_3) \cdot z_2 + (z_2 \cdot z_3) \cdot z_1 - b(z_1, z_2)\, z_3 & (z_1 \cdot z_4) \cdot z_2 + (z_2 \cdot z_4) \cdot z_1 - b(z_1, z_2)\, z_4 \\ (z_2 \cdot z_3) \cdot z_2 - q(z_2)\, z_3 & (z_2 \cdot z_4) \cdot z_2 - q(z_2)\, z_4 \end{pmatrix}$$

# ■ Definitions and rules for proving - first list of identities

```
h[x_] := b[x, x·x]

rules1 = {
   b[x_, y_] b[z_, u_] → b[x·z, y·u] + b[x·u, y·z], (*3*)
   b[x_·y_, x_·z_] → q[x] b[y, z], (*4*)
   b[x_·y_, z_·y_] → q[y] b[x, z], (*5*)
   (x_·y_)·x_ → q[x] y, x_·(y_·x_) → q[x] y, (*6*)
   b[x_*q[y_], z_] → q[y] b[x, z],
   b[x_, z_*q[y_]] → q[y] b[x, z]};

assleft = b[x_, y_·z_] → b[x·y, z];
```

■ $b(s, t^2) = b(t, s^2) = a - c$

```
b[(x·y), (y·x)·(y·x)] -
   b[x·y, y] b[y·x, x] + b[x, y] q[x] q[y] //. rules1
```

$b(x·y, (y·x)·(y·x)) - b((x·y)·(y·x), y·x) - b(y, y·(y·x)) q(x) + b(x, y) q(x) q(y)$

```
% //. assleft
```

$b(x, y) q(x) q(y) - b((y·y)·y, x) q(x)$

```
% //. rules1
```

$b(x, y) q(x) q(y) - b(y, x) q(x) q(y)$

```
% /. b[y, x] → b[x, y]
```

0

■ $b(x s, y t) = e - c$

```
b[x·(x·y), y·(y·x)] -
   b[x, y] b[x·y, y·x] + b[x, y] q[x] q[y] //. rules1
```

$b(x, y) q(x) q(y) - b(y, x) q(x) q(y)$

```
% /. b[y, x] → b[x, y]
```

0

■ $h(s) = a - b(s y, x s)$ (and $h(t) = a - b(t x, y t)$ interchanching $x$ and $y$.

a little change

```
h[x·y] - b[x·y, y] b[x·y, x] + b[(x·y)·y, x·(x·y)] //. rules1
```

$b(x·y, (x·y)·(x·y)) - b((x·y)·(x·y), x·y)$

```
% //. assleft
```

0

- $b(x, y)^3 = b(x y^2, y x^2) + c + e$.

  ```
  b[x, y] b[x · x, y · y] +
    b[x, y] b[x · y, y · x] - b[x · (y · y), y · (x · x)] -
    b[x, y] q[x] q[y] - b[x, y] b[x · y, y · x] //. rules1
  ```

  0

## ■ Main identity

$$N\,(\lambda \oplus x)\,N\,(\mu \oplus y) - N\,((\lambda \oplus x) \bullet (\mu \oplus y)) =$$
$$(1 - \alpha\beta)\,\{3\,b\,(x \cdot (x \cdot y),\,[x,\,y^2]) + (1 + \beta)\,h\,([x,\,y])\} +$$
$$+ 3\,(1 - \alpha\beta)\,b\,([x,\,y],\,-\lambda\,(x \cdot y) \cdot y + \mu\,(y \cdot x) \cdot x + \lambda\mu\,x \cdot y)$$

$$(\lambda \oplus x) \bullet (\mu \oplus y) := (\lambda\,\mu + b(x, y)) \oplus (\lambda\,y + \mu\,x + \alpha\,x \cdot y + \beta\,y \cdot x)$$

```
con[x_, y_] := x · y - y · x;

rules2 = {
   b[x_, y_] b[z_, u_] → b[x · z, y · u] + b[x · u, y · z], (*3*)
   b[x_ · y_, x_ · z_] → q[x] b[y, z], (*4*)
   b[x_ · y_, z_ · y_] -> q[y] b[x, z], (*5*)
   (x_ · y_) · x_ → q[x] y, x_ · (y_ · x_) → q[x] y, (*6*)
   b[x_ * q[y_], z_] → q[y] b[x, z],
   b[x_, z_ * q[y_]] → q[y] b[x, z],
   b[x_, y_]^2 → b[x · x, y · y] + b[x · y, y · x], (*3*)
   b[x_, y_]^3 → b[x, y] b[x · x, y · y] + b[x, y] b[x · y, y · x],
   (*b^3=b b^2*)
   b[x_, x_] → 2 q[x],
   q[x_ · y_] → q[x] q[y] (*1*)};

assocb = {
   b[x_ · y_, (x_ · y_) · (y_ · x_)] → b[(x · y) · (x · y), y · x],
   b[x_ · y_, (y_ · x_) · (x_ · y_)] → b[(x · y) · (x · y), y · x],
   b[y_ · x_, (x_ · y_) · (x_ · y_)] → b[(x · y) · (x · y), y · x]};

leftside = (λ^3 - 3 λ q[x] + h[x]) (μ^3 - 3 μ q[y] + h[y]) -
    ((λ μ + b[x, y])^3 - 3 (λ μ + b[x, y]) q[λ y + μ x + α x · y + β y · x] +
      h[λ y + μ x + α x · y + β y · x]);

rightside =
   (1 - α β) (3 b[x · (x · y), con[x, y · y]] + (1 + β) h[con[x, y]]) +
    3 (1 - α β) b[con[x, y], -λ (x · y) · y + μ (y · x) · x + λ μ x · y];
```

```
dif = Expand[leftside - rightside];
```

$b(y, y \cdot y) \lambda^3 - 3 \mu q(y) \lambda^3 - 3 \mu^2 b(x, y) \lambda^2 - 3 \mu b(x, y)^2 \lambda - 3 \mu^3 q(x) \lambda -$
$\quad 3 b(y, y \cdot y) q(x) \lambda + 9 \mu q(x) q(y) \lambda + 3 \mu q(y \lambda + x \mu + \alpha x \cdot y + \beta y \cdot x) \lambda - b(x, y)^3 +$
$\quad \mu^3 b(x, x \cdot x) + b(x, x \cdot x) b(y, y \cdot y) + 3 \alpha \beta b(x \cdot (x \cdot y), x \cdot (y \cdot y) - (y \cdot y) \cdot x) -$
$\quad 3 b(x \cdot (x \cdot y), x \cdot (y \cdot y) - (y \cdot y) \cdot x) + \alpha \beta^2 b(x \cdot y - y \cdot x, (x \cdot y - y \cdot x) \cdot (x \cdot y - y \cdot x)) +$
$\quad \alpha \beta b(x \cdot y - y \cdot x, (x \cdot y - y \cdot x) \cdot (x \cdot y - y \cdot x)) -$
$\quad \beta b(x \cdot y - y \cdot x, (x \cdot y - y \cdot x) \cdot (x \cdot y - y \cdot x)) - b(x \cdot y - y \cdot x, (x \cdot y - y \cdot x) \cdot (x \cdot y - y \cdot x)) +$
$\quad 3 \alpha \beta b(x \cdot y - y \cdot x, \lambda \mu x \cdot y - \lambda (x \cdot y) \cdot y + \mu (y \cdot x) \cdot x) -$
$\quad 3 b(x \cdot y - y \cdot x, \lambda \mu x \cdot y - \lambda (x \cdot y) \cdot y + \mu (y \cdot x) \cdot x) -$
$\quad b(y \lambda + x \mu + \alpha x \cdot y + \beta y \cdot x, (y \lambda + x \mu + \alpha x \cdot y + \beta y \cdot x) \cdot (y \lambda + x \mu + \alpha x \cdot y + \beta y \cdot x)) -$
$\quad 3 \mu b(x, x \cdot x) q(y) + 3 b(x, y) q(y \lambda + x \mu + \alpha x \cdot y + \beta y \cdot x)$

```
dif //. expandq
```

$b(y, y \cdot y) \lambda^3 - 3 \mu q(y) \lambda^3 - 3 \mu^2 b(x, y) \lambda^2 -$
$\quad 3 \mu b(x, y)^2 \lambda - 3 \mu^3 q(x) \lambda - 3 b(y, y \cdot y) q(x) \lambda + 9 \mu q(x) q(y) \lambda +$
$\quad 3 \mu (b(y \lambda, x \mu + \alpha x \cdot y + \beta y \cdot x) + b(x \mu, \alpha x \cdot y + \beta y \cdot x) +$
$\quad\quad b(\alpha x \cdot y, \beta y \cdot x) + q(y \lambda) + q(x \mu) + q(\alpha x \cdot y) + q(\beta y \cdot x)) \lambda - b(x, y)^3 +$
$\quad \mu^3 b(x, x \cdot x) + b(x, x \cdot x) b(y, y \cdot y) + 3 \alpha \beta b(x \cdot (x \cdot y), x \cdot (y \cdot y) - (y \cdot y) \cdot x) -$
$\quad 3 b(x \cdot (x \cdot y), x \cdot (y \cdot y) - (y \cdot y) \cdot x) + \alpha \beta^2 b(x \cdot y - y \cdot x, (x \cdot y - y \cdot x) \cdot (x \cdot y - y \cdot x)) +$
$\quad \alpha \beta b(x \cdot y - y \cdot x, (x \cdot y - y \cdot x) \cdot (x \cdot y - y \cdot x)) -$
$\quad \beta b(x \cdot y - y \cdot x, (x \cdot y - y \cdot x) \cdot (x \cdot y - y \cdot x)) - b(x \cdot y - y \cdot x, (x \cdot y - y \cdot x) \cdot (x \cdot y - y \cdot x)) +$
$\quad 3 \alpha \beta b(x \cdot y - y \cdot x, \lambda \mu x \cdot y - \lambda (x \cdot y) \cdot y + \mu (y \cdot x) \cdot x) -$
$\quad 3 b(x \cdot y - y \cdot x, \lambda \mu x \cdot y - \lambda (x \cdot y) \cdot y + \mu (y \cdot x) \cdot x) -$
$\quad b(y \lambda + x \mu + \alpha x \cdot y + \beta y \cdot x, (y \lambda + x \mu + \alpha x \cdot y + \beta y \cdot x) \cdot (y \lambda + x \mu + \alpha x \cdot y + \beta y \cdot x)) -$
$\quad 3 \mu b(x, x \cdot x) q(y) + 3 b(x, y) (b(y \lambda, x \mu + \alpha x \cdot y + \beta y \cdot x) + b(x \mu, \alpha x \cdot y + \beta y \cdot x) +$
$\quad\quad b(\alpha x \cdot y, \beta y \cdot x) + q(y \lambda) + q(x \mu) + q(\alpha x \cdot y) + q(\beta y \cdot x))$

```
% //. expanddot
```

$b(y, y \cdot y) \lambda^3 - 3 \mu q(y) \lambda^3 - 3 \mu^2 b(x, y) \lambda^2 -$
$\quad 3 \mu b(x, y)^2 \lambda - 3 \mu^3 q(x) \lambda - 3 b(y, y \cdot y) q(x) \lambda + 9 \mu q(x) q(y) \lambda +$
$\quad 3 \mu (b(y \lambda, x \mu + \alpha x \cdot y + \beta y \cdot x) + b(x \mu, \alpha x \cdot y + \beta y \cdot x) +$
$\quad\quad b(\alpha x \cdot y, \beta y \cdot x) + q(y \lambda) + q(x \mu) + q(\alpha x \cdot y) + q(\beta y \cdot x)) \lambda - b(x, y)^3 +$
$\quad \mu^3 b(x, x \cdot x) + b(x, x \cdot x) b(y, y \cdot y) + 3 \alpha \beta b(x \cdot (x \cdot y), x \cdot (y \cdot y) - (y \cdot y) \cdot x) -$
$\quad 3 b(x \cdot (x \cdot y), x \cdot (y \cdot y) - (y \cdot y) \cdot x) +$
$\quad 3 \alpha \beta b(x \cdot y - y \cdot x, \lambda \mu x \cdot y - \lambda (x \cdot y) \cdot y + \mu (y \cdot x) \cdot x) -$
$\quad 3 b(x \cdot y - y \cdot x, \lambda \mu x \cdot y - \lambda (x \cdot y) \cdot y + \mu (y \cdot x) \cdot x) +$
$\quad \alpha \beta^2 b(x \cdot y - y \cdot x, -(x \cdot y) \cdot (y \cdot x) + (x \cdot y) \cdot (x \cdot y) + -(y \cdot x) \cdot (x \cdot y) + (y \cdot x) \cdot (y \cdot x)) +$
$\quad \alpha \beta b(x \cdot y - y \cdot x, -(x \cdot y) \cdot (y \cdot x) + (x \cdot y) \cdot (x \cdot y) + -(y \cdot x) \cdot (x \cdot y) + (y \cdot x) \cdot (y \cdot x)) -$
$\quad \beta b(x \cdot y - y \cdot x, -(x \cdot y) \cdot (y \cdot x) + (x \cdot y) \cdot (x \cdot y) + -(y \cdot x) \cdot (x \cdot y) + (y \cdot x) \cdot (y \cdot x)) -$
$\quad b(x \cdot y - y \cdot x, -(x \cdot y) \cdot (y \cdot x) + (x \cdot y) \cdot (x \cdot y) + -(y \cdot x) \cdot (x \cdot y) + (y \cdot x) \cdot (y \cdot x)) -$
$\quad b(y \lambda + x \mu + \alpha x \cdot y + \beta y \cdot x, (y \lambda) \cdot (y \lambda) + (y \lambda) \cdot (x \mu) + (y \lambda) \cdot (\alpha x \cdot y) +$
$\quad\quad (y \lambda) \cdot (\beta y \cdot x) + (x \mu) \cdot (y \lambda) + (x \mu) \cdot (x \mu) + (x \mu) \cdot (\alpha x \cdot y) + (x \mu) \cdot (\beta y \cdot x) +$
$\quad\quad (\alpha x \cdot y) \cdot (y \lambda) + (\alpha x \cdot y) \cdot (x \mu) + (\alpha x \cdot y) \cdot (\alpha x \cdot y) + (\alpha x \cdot y) \cdot (\beta y \cdot x) +$
$\quad\quad (\beta y \cdot x) \cdot (y \lambda) + (\beta y \cdot x) \cdot (x \mu) + (\beta y \cdot x) \cdot (\alpha x \cdot y) + (\beta y \cdot x) \cdot (\beta y \cdot x)) -$
$\quad 3 \mu b(x, x \cdot x) q(y) + 3 b(x, y) (b(y \lambda, x \mu + \alpha x \cdot y + \beta y \cdot x) + b(x \mu, \alpha x \cdot y + \beta y \cdot x) +$
$\quad\quad b(\alpha x \cdot y, \beta y \cdot x) + q(y \lambda) + q(x \mu) + q(\alpha x \cdot y) + q(\beta y \cdot x))$

**% //. expandb**

$b(y, y \cdot y) \lambda^3 - 3 \mu q(y) \lambda^3 - 3 \mu^2 b(x, y) \lambda^2 - 3 \mu b(x, y)^2 \lambda - 3 \mu^3 q(x) \lambda - 3 b(y, y \cdot y) q(x) \lambda + 9 \mu q(x) q(y) \lambda + 3 \mu (b(y \lambda, x \mu) + b(y \lambda, \alpha x \cdot y) + b(y \lambda, \beta y \cdot x) + b(x \mu, \alpha x \cdot y) + b(x \mu, \beta y \cdot x) + b(\alpha x \cdot y, \beta y \cdot x) + q(y \lambda) + q(x \mu) + q(\alpha x \cdot y) + q(\beta y \cdot x)) \lambda - b(x, y)^3 + \mu^3 b(x, x \cdot x) + b(x, x \cdot x) b(y, y \cdot y) - b(y \lambda, (y \lambda) \cdot (y \lambda)) - b(y \lambda, (y \lambda) \cdot (x \mu)) - b(y \lambda, (y \lambda) \cdot (\alpha x \cdot y)) - b(y \lambda, (y \lambda) \cdot (\beta y \cdot x)) - b(y \lambda, (x \mu) \cdot (y \lambda)) - b(y \lambda, (x \mu) \cdot (x \mu)) - b(y \lambda, (x \mu) \cdot (\alpha x \cdot y)) - b(y \lambda, (x \mu) \cdot (\beta y \cdot x)) - b(y \lambda, (\alpha x \cdot y) \cdot (y \lambda)) - b(y \lambda, (\alpha x \cdot y) \cdot (x \mu)) - b(y \lambda, (\alpha x \cdot y) \cdot (\alpha x \cdot y)) - b(y \lambda, (\alpha x \cdot y) \cdot (\beta y \cdot x)) - b(y \lambda, (\beta y \cdot x) \cdot (y \lambda)) - b(y \lambda, (\beta y \cdot x) \cdot (x \mu)) - b(y \lambda, (\beta y \cdot x) \cdot (\alpha x \cdot y)) - b(y \lambda, (\beta y \cdot x) \cdot (\beta y \cdot x)) - b(x \mu, (y \lambda) \cdot (y \lambda)) - b(x \mu, (y \lambda) \cdot (x \mu)) - b(x \mu, (y \lambda) \cdot (\alpha x \cdot y)) - b(x \mu, (y \lambda) \cdot (\beta y \cdot x)) - b(x \mu, (x \mu) \cdot (y \lambda)) - b(x \mu, (x \mu) \cdot (x \mu)) - b(x \mu, (x \mu) \cdot (\alpha x \cdot y)) - b(x \mu, (x \mu) \cdot (\beta y \cdot x)) - b(x \mu, (\alpha x \cdot y) \cdot (y \lambda)) - b(x \mu, (\alpha x \cdot y) \cdot (x \mu)) - b(x \mu, (\alpha x \cdot y) \cdot (\alpha x \cdot y)) - b(x \mu, (\alpha x \cdot y) \cdot (\beta y \cdot x)) - b(x \mu, (\beta y \cdot x) \cdot (y \lambda)) - b(x \mu, (\beta y \cdot x) \cdot (x \mu)) - b(x \mu, (\beta y \cdot x) \cdot (\alpha x \cdot y)) - b(x \mu, (\beta y \cdot x) \cdot (\beta y \cdot x)) - b(x \cdot y, (x \cdot y) \cdot (x \cdot y)) + b(x \cdot y, (x \cdot y) \cdot (y \cdot x)) + b(x \cdot y, (y \cdot x) \cdot (x \cdot y)) - b(x \cdot y, (y \cdot x) \cdot (y \cdot x)) - b(\alpha x \cdot y, (y \lambda) \cdot (y \lambda)) - b(\alpha x \cdot y, (y \lambda) \cdot (x \mu)) - b(\alpha x \cdot y, (y \lambda) \cdot (\alpha x \cdot y)) - b(\alpha x \cdot y, (y \lambda) \cdot (\beta y \cdot x)) - b(\alpha x \cdot y, (x \mu) \cdot (y \lambda)) - b(\alpha x \cdot y, (x \mu) \cdot (x \mu)) - b(\alpha x \cdot y, (x \mu) \cdot (\alpha x \cdot y)) - b(\alpha x \cdot y, (x \mu) \cdot (\beta y \cdot x)) - b(\alpha x \cdot y, (\alpha x \cdot y) \cdot (y \lambda)) - b(\alpha x \cdot y, (\alpha x \cdot y) \cdot (x \mu)) - b(\alpha x \cdot y, (\alpha x \cdot y) \cdot (\alpha x \cdot y)) - b(\alpha x \cdot y, (\alpha x \cdot y) \cdot (\beta y \cdot x)) - b(\alpha x \cdot y, (\beta y \cdot x) \cdot (y \lambda)) - b(\alpha x \cdot y, (\beta y \cdot x) \cdot (x \mu)) - b(\alpha x \cdot y, (\beta y \cdot x) \cdot (\alpha x \cdot y)) - b(\alpha x \cdot y, (\beta y \cdot x) \cdot (\beta y \cdot x)) + 3 \alpha \beta (b(x \cdot (x \cdot y), x \cdot (y \cdot y)) - b(x \cdot (x \cdot y), (y \cdot y) \cdot x)) - 3 (b(x \cdot (x \cdot y), x \cdot (y \cdot y)) - b(x \cdot (x \cdot y), (y \cdot y) \cdot x)) + b(y \cdot x, (x \cdot y) \cdot (x \cdot y)) - b(y \cdot x, (x \cdot y) \cdot (y \cdot x)) + 3 \alpha \beta (b(x \cdot y, \lambda \mu x \cdot y) + b(x \cdot y, -\lambda (x \cdot y) \cdot y) + b(x \cdot y, \mu (y \cdot x) \cdot x) - b(y \cdot x, \lambda \mu x \cdot y) - b(y \cdot x, -\lambda (x \cdot y) \cdot y) - b(y \cdot x, \mu (y \cdot x) \cdot x)) - 3 (b(x \cdot y, \lambda \mu x \cdot y) + b(x \cdot y, -\lambda (x \cdot y) \cdot y) + b(x \cdot y, \mu (y \cdot x) \cdot x) - b(y \cdot x, \lambda \mu x \cdot y) - b(y \cdot x, -\lambda (x \cdot y) \cdot y) - b(y \cdot x, \mu (y \cdot x) \cdot x)) - b(y \cdot x, (y \cdot x) \cdot (x \cdot y)) + \alpha \beta^2 (b(x \cdot y, (x \cdot y) \cdot (x \cdot y)) - b(x \cdot y, (x \cdot y) \cdot (y \cdot x)) - b(x \cdot y, (y \cdot x) \cdot (x \cdot y)) + b(x \cdot y, (y \cdot x) \cdot (y \cdot x)) - b(y \cdot x, (x \cdot y) \cdot (x \cdot y)) + b(y \cdot x, (x \cdot y) \cdot (y \cdot x)) + b(y \cdot x, (y \cdot x) \cdot (x \cdot y)) - b(y \cdot x, (y \cdot x) \cdot (y \cdot x))) + \alpha \beta (b(x \cdot y, (x \cdot y) \cdot (x \cdot y)) - b(x \cdot y, (x \cdot y) \cdot (y \cdot x)) - b(x \cdot y, (y \cdot x) \cdot (x \cdot y)) + b(x \cdot y, (y \cdot x) \cdot (y \cdot x)) - b(y \cdot x, (x \cdot y) \cdot (x \cdot y)) + b(y \cdot x, (x \cdot y) \cdot (y \cdot x)) + b(y \cdot x, (y \cdot x) \cdot (x \cdot y)) - b(y \cdot x, (y \cdot x) \cdot (y \cdot x))) - \beta (b(x \cdot y, (x \cdot y) \cdot (x \cdot y)) - b(x \cdot y, (x \cdot y) \cdot (y \cdot x)) - b(x \cdot y, (y \cdot x) \cdot (x \cdot y)) + b(x \cdot y, (y \cdot x) \cdot (y \cdot x)) - b(y \cdot x, (x \cdot y) \cdot (x \cdot y)) + b(y \cdot x, (x \cdot y) \cdot (y \cdot x)) + b(y \cdot x, (y \cdot x) \cdot (x \cdot y)) - b(y \cdot x, (y \cdot x) \cdot (y \cdot x))) + b(y \cdot x, (y \cdot x) \cdot (y \cdot x)) - b(\beta y \cdot x, (y \lambda) \cdot (y \lambda)) - b(\beta y \cdot x, (y \lambda) \cdot (x \mu)) - b(\beta y \cdot x, (y \lambda) \cdot (\alpha x \cdot y)) - b(\beta y \cdot x, (y \lambda) \cdot (\beta y \cdot x)) - b(\beta y \cdot x, (x \mu) \cdot (y \lambda)) - b(\beta y \cdot x, (x \mu) \cdot (x \mu)) - b(\beta y \cdot x, (x \mu) \cdot (\alpha x \cdot y)) - b(\beta y \cdot x, (x \mu) \cdot (\beta y \cdot x)) - b(\beta y \cdot x, (\alpha x \cdot y) \cdot (y \lambda)) - b(\beta y \cdot x, (\alpha x \cdot y) \cdot (x \mu)) - b(\beta y \cdot x, (\alpha x \cdot y) \cdot (\alpha x \cdot y)) - b(\beta y \cdot x, (\alpha x \cdot y) \cdot (\beta y \cdot x)) - b(\beta y \cdot x, (\beta y \cdot x) \cdot (y \lambda)) - b(\beta y \cdot x, (\beta y \cdot x) \cdot (x \mu)) - b(\beta y \cdot x, (\beta y \cdot x) \cdot (\alpha x \cdot y)) - b(\beta y \cdot x, (\beta y \cdot x) \cdot (\beta y \cdot x)) - 3 \mu b(x, x \cdot x) q(y) + 3 b(x, y) (b(y \lambda, x \mu) + b(y \lambda, \alpha x \cdot y) + b(y \lambda, \beta y \cdot x) + b(x \mu, \alpha x \cdot y) + b(x \mu, \beta y \cdot x) + b(\alpha x \cdot y, \beta y \cdot x) + q(y \lambda) + q(x \mu) + q(\alpha x \cdot y) + q(\beta y \cdot x))$

**Expand[%]**

$b(y, y \cdot y) \lambda^3 - 3 \mu q(y) \lambda^3 - 3 \mu^2 b(x, y) \lambda^2 - 3 \mu b(x, y)^2 \lambda + 3 \mu b(y \lambda, x \mu) \lambda +$
$3 \mu b(y \lambda, \alpha x \cdot y) \lambda + 3 \mu b(y \lambda, \beta y \cdot x) \lambda + 3 \mu b(x \mu, \alpha x \cdot y) \lambda + 3 \mu b(x \mu, \beta y \cdot x) \lambda +$
$3 \mu b(\alpha x \cdot y, \beta y \cdot x) \lambda - 3 \mu^3 q(x) \lambda - 3 b(y, y \cdot y) q(x) \lambda + 9 \mu q(x) q(y) \lambda +$
$3 \mu q(y \lambda) \lambda + 3 \mu q(x \mu) \lambda + 3 \mu q(\alpha x \cdot y) \lambda + 3 \mu q(\beta y \cdot x) \lambda - b(x, y)^3 +$
$\mu^3 b(x, x \cdot x) + b(x, x \cdot x) b(y, y \cdot y) + 3 b(x, y) b(y \lambda, x \mu) + 3 b(x, y) b(y \lambda, \alpha x \cdot y) +$
$3 b(x, y) b(y \lambda, \beta y \cdot x) - b(y \lambda, (y \lambda) \cdot (y \lambda)) - b(y \lambda, (y \lambda) \cdot (x \mu)) -$
$b(y \lambda, (y \lambda) \cdot (\alpha x \cdot y)) - b(y \lambda, (y \lambda) \cdot (\beta y \cdot x)) - b(y \lambda, (x \mu) \cdot (y \lambda)) - b(y \lambda, (x \mu) \cdot (x \mu)) -$
$b(y \lambda, (x \mu) \cdot (\alpha x \cdot y)) - b(y \lambda, (x \mu) \cdot (\beta y \cdot x)) - b(y \lambda, (\alpha x \cdot y) \cdot (y \lambda)) -$
$b(y \lambda, (\alpha x \cdot y) \cdot (x \mu)) - b(y \lambda, (\alpha x \cdot y) \cdot (\alpha x \cdot y)) - b(y \lambda, (\alpha x \cdot y) \cdot (\beta y \cdot x)) -$
$b(y \lambda, (\beta y \cdot x) \cdot (y \lambda)) - b(y \lambda, (\beta y \cdot x) \cdot (x \mu)) - b(y \lambda, (\beta y \cdot x) \cdot (\alpha x \cdot y)) -$
$b(y \lambda, (\beta y \cdot x) \cdot (\beta y \cdot x)) + 3 b(x, y) b(x \mu, \alpha x \cdot y) + 3 b(x, y) b(x \mu, \beta y \cdot x) -$
$b(x \mu, (y \lambda) \cdot (y \lambda)) - b(x \mu, (y \lambda) \cdot (x \mu)) - b(x \mu, (y \lambda) \cdot (\alpha x \cdot y)) - b(x \mu, (y \lambda) \cdot (\beta y \cdot x)) -$
$b(x \mu, (x \mu) \cdot (y \lambda)) - b(x \mu, (x \mu) \cdot (x \mu)) - b(x \mu, (x \mu) \cdot (\alpha x \cdot y)) - b(x \mu, (x \mu) \cdot (\beta y \cdot x)) -$
$b(x \mu, (\alpha x \cdot y) \cdot (y \lambda)) - b(x \mu, (\alpha x \cdot y) \cdot (x \mu)) - b(x \mu, (\alpha x \cdot y) \cdot (\alpha x \cdot y)) -$
$b(x \mu, (\alpha x \cdot y) \cdot (\beta y \cdot x)) - b(x \mu, (\beta y \cdot x) \cdot (y \lambda)) - b(x \mu, (\beta y \cdot x) \cdot (x \mu)) -$
$b(x \mu, (\beta y \cdot x) \cdot (\alpha x \cdot y)) - b(x \mu, (\beta y \cdot x) \cdot (\beta y \cdot x)) + 3 \alpha \beta b(x \cdot y, \lambda \mu x \cdot y) -$
$3 b(x \cdot y, \lambda \mu x \cdot y) + 3 \alpha \beta b(x \cdot y, -\lambda (x \cdot y) \cdot y) - 3 b(x \cdot y, -\lambda (x \cdot y) \cdot y) +$
$\alpha \beta^2 b(x \cdot y, (x \cdot y) \cdot (x \cdot y)) + \alpha \beta b(x \cdot y, (x \cdot y) \cdot (x \cdot y)) - \beta b(x \cdot y, (x \cdot y) \cdot (x \cdot y)) -$
$b(x \cdot y, (x \cdot y) \cdot (x \cdot y)) - \alpha \beta^2 b(x \cdot y, (x \cdot y) \cdot (y \cdot x)) - \alpha \beta b(x \cdot y, (x \cdot y) \cdot (y \cdot x)) +$
$\beta b(x \cdot y, (x \cdot y) \cdot (y \cdot x)) + b(x \cdot y, (x \cdot y) \cdot (y \cdot x)) + 3 \alpha \beta b(x \cdot y, \mu (y \cdot x) \cdot x) -$
$3 b(x \cdot y, \mu (y \cdot x) \cdot x) - \alpha \beta^2 b(x \cdot y, (y \cdot x) \cdot (x \cdot y)) - \alpha \beta b(x \cdot y, (y \cdot x) \cdot (x \cdot y)) +$
$\beta b(x \cdot y, (y \cdot x) \cdot (x \cdot y)) + b(x \cdot y, (y \cdot x) \cdot (x \cdot y)) + \alpha \beta^2 b(x \cdot y, (y \cdot x) \cdot (y \cdot x)) +$
$\alpha \beta b(x \cdot y, (y \cdot x) \cdot (y \cdot x)) - \beta b(x \cdot y, (y \cdot x) \cdot (y \cdot x)) - b(x \cdot y, (y \cdot x) \cdot (y \cdot x)) +$
$3 b(x, y) b(\alpha x \cdot y, \beta y \cdot x) - b(\alpha x \cdot y, (y \lambda) \cdot (y \lambda)) - b(\alpha x \cdot y, (y \lambda) \cdot (x \mu)) -$
$b(\alpha x \cdot y, (y \lambda) \cdot (\alpha x \cdot y)) - b(\alpha x \cdot y, (y \lambda) \cdot (\beta y \cdot x)) - b(\alpha x \cdot y, (x \mu) \cdot (y \lambda)) -$
$b(\alpha x \cdot y, (x \mu) \cdot (x \mu)) - b(\alpha x \cdot y, (x \mu) \cdot (\alpha x \cdot y)) - b(\alpha x \cdot y, (x \mu) \cdot (\beta y \cdot x)) -$
$b(\alpha x \cdot y, (\alpha x \cdot y) \cdot (y \lambda)) - b(\alpha x \cdot y, (\alpha x \cdot y) \cdot (x \mu)) - b(\alpha x \cdot y, (\alpha x \cdot y) \cdot (\alpha x \cdot y)) -$
$b(\alpha x \cdot y, (\alpha x \cdot y) \cdot (\beta y \cdot x)) - b(\alpha x \cdot y, (\beta y \cdot x) \cdot (y \lambda)) - b(\alpha x \cdot y, (\beta y \cdot x) \cdot (x \mu)) -$
$b(\alpha x \cdot y, (\beta y \cdot x) \cdot (\alpha x \cdot y)) - b(\alpha x \cdot y, (\beta y \cdot x) \cdot (\beta y \cdot x)) + 3 \alpha \beta b(x \cdot (x \cdot y), x \cdot (y \cdot y)) -$
$3 b(x \cdot (x \cdot y), x \cdot (y \cdot y)) - 3 \alpha \beta b(x \cdot (x \cdot y), (y \cdot y) \cdot x) + 3 b(x \cdot (x \cdot y), (y \cdot y) \cdot x) -$
$3 \alpha \beta b(y \cdot x, \lambda \mu x \cdot y) + 3 b(y \cdot x, \lambda \mu x \cdot y) - 3 \alpha \beta b(y \cdot x, -\lambda (x \cdot y) \cdot y) +$
$3 b(y \cdot x, -\lambda (x \cdot y) \cdot y) - \alpha \beta^2 b(y \cdot x, (x \cdot y) \cdot (x \cdot y)) - \alpha \beta b(y \cdot x, (x \cdot y) \cdot (x \cdot y)) +$
$\beta b(y \cdot x, (x \cdot y) \cdot (x \cdot y)) + b(y \cdot x, (x \cdot y) \cdot (x \cdot y)) + \alpha \beta^2 b(y \cdot x, (x \cdot y) \cdot (y \cdot x)) +$
$\alpha \beta b(y \cdot x, (x \cdot y) \cdot (y \cdot x)) - \beta b(y \cdot x, (x \cdot y) \cdot (y \cdot x)) - b(y \cdot x, (x \cdot y) \cdot (y \cdot x)) -$
$3 \alpha \beta b(y \cdot x, \mu (y \cdot x) \cdot x) + 3 b(y \cdot x, \mu (y \cdot x) \cdot x) + \alpha \beta^2 b(y \cdot x, (y \cdot x) \cdot (x \cdot y)) +$
$\alpha \beta b(y \cdot x, (y \cdot x) \cdot (x \cdot y)) - \beta b(y \cdot x, (y \cdot x) \cdot (x \cdot y)) - b(y \cdot x, (y \cdot x) \cdot (x \cdot y)) -$
$\alpha \beta^2 b(y \cdot x, (y \cdot x) \cdot (y \cdot x)) - \alpha \beta b(y \cdot x, (y \cdot x) \cdot (y \cdot x)) + \beta b(y \cdot x, (y \cdot x) \cdot (y \cdot x)) +$
$b(y \cdot x, (y \cdot x) \cdot (y \cdot x)) - b(\beta y \cdot x, (y \lambda) \cdot (y \lambda)) - b(\beta y \cdot x, (y \lambda) \cdot (x \mu)) -$
$b(\beta y \cdot x, (y \lambda) \cdot (\alpha x \cdot y)) - b(\beta y \cdot x, (y \lambda) \cdot (\beta y \cdot x)) - b(\beta y \cdot x, (x \mu) \cdot (y \lambda)) -$
$b(\beta y \cdot x, (x \mu) \cdot (x \mu)) - b(\beta y \cdot x, (x \mu) \cdot (\alpha x \cdot y)) - b(\beta y \cdot x, (x \mu) \cdot (\beta y \cdot x)) -$
$b(\beta y \cdot x, (\alpha x \cdot y) \cdot (y \lambda)) - b(\beta y \cdot x, (\alpha x \cdot y) \cdot (x \mu)) - b(\beta y \cdot x, (\alpha x \cdot y) \cdot (\alpha x \cdot y)) -$
$b(\beta y \cdot x, (\alpha x \cdot y) \cdot (\beta y \cdot x)) - b(\beta y \cdot x, (\beta y \cdot x) \cdot (y \lambda)) - b(\beta y \cdot x, (\beta y \cdot x) \cdot (x \mu)) -$
$b(\beta y \cdot x, (\beta y \cdot x) \cdot (\alpha x \cdot y)) - b(\beta y \cdot x, (\beta y \cdot x) \cdot (\beta y \cdot x)) - 3 \mu b(x, x \cdot x) q(y) +$
$3 b(x, y) q(y \lambda) + 3 b(x, y) q(x \mu) + 3 b(x, y) q(\alpha x \cdot y) + 3 b(x, y) q(\beta y \cdot x)$

**% //. rules2**

$b(y, y \cdot y) \lambda^3 - 3 \mu q(y) \lambda^3 - 3 \mu^2 b(x, y) \lambda^2 + 3 \mu b(y\lambda, x\mu) \lambda +$
$3 \mu b(y\lambda, \alpha x \cdot y) \lambda + 3 \mu b(y\lambda, \beta y \cdot x) \lambda + 3 \mu b(x\mu, \alpha x \cdot y) \lambda + 3 \mu b(x\mu, \beta y \cdot x) \lambda -$
$3 \mu (b(x \cdot x, y \cdot y) + b(x \cdot y, y \cdot x)) \lambda + 3 \mu b(\alpha x \cdot y, \beta y \cdot x) \lambda - 3 \mu^3 q(x) \lambda -$
$3 b(y, y \cdot y) q(x) \lambda + 9 \mu q(x) q(y) \lambda + 3 \mu q(y\lambda) \lambda + 3 \mu q(x\mu) \lambda + 3 \mu q(\alpha x \cdot y) \lambda +$
$3 \mu q(\beta y \cdot x) \lambda + \mu^3 b(x, x \cdot x) - b(y\lambda, (y\lambda) \cdot (y\lambda)) - b(y\lambda, (y\lambda) \cdot (x\mu)) -$
$b(y\lambda, (y\lambda) \cdot (\alpha x \cdot y)) - b(y\lambda, (y\lambda) \cdot (\beta y \cdot x)) - b(y\lambda, (x\mu) \cdot (y\lambda)) - b(y\lambda, (x\mu) \cdot (x\mu)) -$
$b(y\lambda, (x\mu) \cdot (\alpha x \cdot y)) - b(y\lambda, (x\mu) \cdot (\beta y \cdot x)) - b(y\lambda, (\alpha x \cdot y) \cdot (y\lambda)) -$
$b(y\lambda, (\alpha x \cdot y) \cdot (x\mu)) - b(y\lambda, (\alpha x \cdot y) \cdot (\alpha x \cdot y)) - b(y\lambda, (\alpha x \cdot y) \cdot (\beta y \cdot x)) -$
$b(y\lambda, (\beta y \cdot x) \cdot (y\lambda)) - b(y\lambda, (\beta y \cdot x) \cdot (x\mu)) - b(y\lambda, (\beta y \cdot x) \cdot (\alpha x \cdot y)) -$
$b(y\lambda, (\beta y \cdot x) \cdot (\beta y \cdot x)) - b(x\mu, (y\lambda) \cdot (y\lambda)) - b(x\mu, (y\lambda) \cdot (x\mu)) -$
$b(x\mu, (y\lambda) \cdot (\alpha x \cdot y)) - b(x\mu, (y\lambda) \cdot (\beta y \cdot x)) - b(x\mu, (x\mu) \cdot (y\lambda)) - b(x\mu, (x\mu) \cdot (x\mu)) -$
$b(x\mu, (x\mu) \cdot (\alpha x \cdot y)) - b(x\mu, (x\mu) \cdot (\beta y \cdot x)) - b(x\mu, (\alpha x \cdot y) \cdot (y\lambda)) -$
$b(x\mu, (\alpha x \cdot y) \cdot (x\mu)) - b(x\mu, (\alpha x \cdot y) \cdot (\alpha x \cdot y)) - b(x\mu, (\alpha x \cdot y) \cdot (\beta y \cdot x)) -$
$b(x\mu, (\beta y \cdot x) \cdot (y\lambda)) - b(x\mu, (\beta y \cdot x) \cdot (x\mu)) - b(x\mu, (\beta y \cdot x) \cdot (\alpha x \cdot y)) -$
$b(x\mu, (\beta y \cdot x) \cdot (\beta y \cdot x)) + 3 \alpha \beta b(x \cdot y, \lambda \mu x \cdot y) - 3 b(x \cdot y, \lambda \mu x \cdot y) +$
$b(x \cdot y, (x \cdot x) \cdot (y \cdot y)) + 3 \alpha \beta b(x \cdot y, -\lambda(x \cdot y) \cdot y) - 3 b(x \cdot y, -\lambda(x \cdot y) \cdot y) +$
$\alpha \beta^2 b(x \cdot y, (x \cdot y) \cdot (x \cdot y)) + \alpha \beta b(x \cdot y, (x \cdot y) \cdot (x \cdot y)) - \beta b(x \cdot y, (x \cdot y) \cdot (x \cdot y)) -$
$b(x \cdot y, (x \cdot y) \cdot (x \cdot y)) - \alpha \beta^2 b(x \cdot y, (x \cdot y) \cdot (y \cdot x)) - \alpha \beta b(x \cdot y, (x \cdot y) \cdot (y \cdot x)) +$
$\beta b(x \cdot y, (x \cdot y) \cdot (y \cdot x)) + b(x \cdot y, (x \cdot y) \cdot (y \cdot x)) + 3 \alpha \beta b(x \cdot y, \mu(y \cdot x) \cdot x) -$
$3 b(x \cdot y, \mu(y \cdot x) \cdot x) - \alpha \beta^2 b(x \cdot y, (y \cdot x) \cdot (x \cdot y)) - \alpha \beta b(x \cdot y, (y \cdot x) \cdot (x \cdot y)) +$
$\beta b(x \cdot y, (y \cdot x) \cdot (x \cdot y)) + b(x \cdot y, (y \cdot x) \cdot (x \cdot y)) + \alpha \beta^2 b(x \cdot y, (y \cdot x) \cdot (y \cdot x)) +$
$\alpha \beta b(x \cdot y, (y \cdot x) \cdot (y \cdot x)) - \beta b(x \cdot y, (y \cdot x) \cdot (y \cdot x)) - b(x \cdot y, (y \cdot x) \cdot (y \cdot x)) -$
$b(\alpha x \cdot y, (y\lambda) \cdot (y\lambda)) - b(\alpha x \cdot y, (y\lambda) \cdot (x\mu)) - b(\alpha x \cdot y, (y\lambda) \cdot (\alpha x \cdot y)) -$
$b(\alpha x \cdot y, (y\lambda) \cdot (\beta y \cdot x)) - b(\alpha x \cdot y, (x\mu) \cdot (y\lambda)) - b(\alpha x \cdot y, (x\mu) \cdot (x\mu)) -$
$b(\alpha x \cdot y, (x\mu) \cdot (\alpha x \cdot y)) - b(\alpha x \cdot y, (x\mu) \cdot (\beta y \cdot x)) - b(\alpha x \cdot y, (\alpha x \cdot y) \cdot (y\lambda)) -$
$b(\alpha x \cdot y, (\alpha x \cdot y) \cdot (x\mu)) - b(\alpha x \cdot y, (\alpha x \cdot y) \cdot (\alpha x \cdot y)) - b(\alpha x \cdot y, (\alpha x \cdot y) \cdot (\beta y \cdot x)) -$
$b(\alpha x \cdot y, (\beta y \cdot x) \cdot (y\lambda)) - b(\alpha x \cdot y, (\beta y \cdot x) \cdot (x\mu)) - b(\alpha x \cdot y, (\beta y \cdot x) \cdot (\alpha x \cdot y)) -$
$b(\alpha x \cdot y, (\beta y \cdot x) \cdot (\beta y \cdot x)) + 3 (b(x \cdot (y\lambda), y \cdot (x\mu)) + b(x \cdot (x\mu), y \cdot (y\lambda))) -$
$b(x \cdot (x \cdot y), y \cdot (y \cdot x)) - 3 \alpha \beta b(x \cdot (x \cdot y), (y \cdot y) \cdot x) + 3 b(x \cdot (x \cdot y), (y \cdot y) \cdot x) +$
$3 (b(x \cdot (y\lambda), y \cdot (\alpha x \cdot y)) + b(x \cdot (\alpha x \cdot y), y \cdot (y\lambda))) +$
$3 (b(x \cdot (x\mu), y \cdot (\alpha x \cdot y)) + b(x \cdot (\alpha x \cdot y), y \cdot (x\mu))) +$
$3 (b(x \cdot (y\lambda), y \cdot (\beta y \cdot x)) + b(x \cdot (\beta y \cdot x), y \cdot (y\lambda))) +$
$3 (b(x \cdot (x\mu), y \cdot (\beta y \cdot x)) + b(x \cdot (\beta y \cdot x), y \cdot (x\mu))) +$
$3 (b(x \cdot (\alpha x \cdot y), y \cdot (\beta y \cdot x)) + b(x \cdot (\beta y \cdot x), y \cdot (\alpha x \cdot y))) - b(x \cdot (y \cdot y), y \cdot (x \cdot x)) +$
$b(x \cdot (y \cdot y), (x \cdot x) \cdot y) - 3 \alpha \beta b(y \cdot x, \lambda \mu x \cdot y) + 3 b(y \cdot x, \lambda \mu x \cdot y) -$
$3 \alpha \beta b(y \cdot x, -\lambda(x \cdot y) \cdot y) + 3 b(y \cdot x, -\lambda(x \cdot y) \cdot y) - \alpha \beta^2 b(y \cdot x, (x \cdot y) \cdot (x \cdot y)) -$
$\alpha \beta b(y \cdot x, (x \cdot y) \cdot (x \cdot y)) + \beta b(y \cdot x, (x \cdot y) \cdot (x \cdot y)) + b(y \cdot x, (x \cdot y) \cdot (x \cdot y)) +$
$\alpha \beta^2 b(y \cdot x, (x \cdot y) \cdot (y \cdot x)) + \alpha \beta b(y \cdot x, (x \cdot y) \cdot (y \cdot x)) - \beta b(y \cdot x, (x \cdot y) \cdot (y \cdot x)) -$
$b(y \cdot x, (x \cdot y) \cdot (y \cdot x)) - 3 \alpha \beta b(y \cdot x, \mu(y \cdot x) \cdot x) + 3 b(y \cdot x, \mu(y \cdot x) \cdot x) +$
$\alpha \beta^2 b(y \cdot x, (y \cdot x) \cdot (x \cdot y)) + \alpha \beta b(y \cdot x, (y \cdot x) \cdot (x \cdot y)) - \beta b(y \cdot x, (y \cdot x) \cdot (x \cdot y)) -$
$b(y \cdot x, (y \cdot x) \cdot (x \cdot y)) - \alpha \beta^2 b(y \cdot x, (y \cdot x) \cdot (y \cdot x)) - \alpha \beta b(y \cdot x, (y \cdot x) \cdot (y \cdot x)) +$
$\beta b(y \cdot x, (y \cdot x) \cdot (y \cdot x)) + b(y \cdot x, (y \cdot x) \cdot (y \cdot x)) - b(\beta y \cdot x, (y\lambda) \cdot (y\lambda)) -$
$b(\beta y \cdot x, (y\lambda) \cdot (x\mu)) - b(\beta y \cdot x, (y\lambda) \cdot (\alpha x \cdot y)) - b(\beta y \cdot x, (y\lambda) \cdot (\beta y \cdot x)) -$
$b(\beta y \cdot x, (x\mu) \cdot (y\lambda)) - b(\beta y \cdot x, (x\mu) \cdot (x\mu)) - b(\beta y \cdot x, (x\mu) \cdot (\alpha x \cdot y)) -$
$b(\beta y \cdot x, (x\mu) \cdot (\beta y \cdot x)) - b(\beta y \cdot x, (\alpha x \cdot y) \cdot (y\lambda)) - b(\beta y \cdot x, (\alpha x \cdot y) \cdot (x\mu)) -$
$b(\beta y \cdot x, (\alpha x \cdot y) \cdot (\alpha x \cdot y)) - b(\beta y \cdot x, (\alpha x \cdot y) \cdot (\beta y \cdot x)) - b(\beta y \cdot x, (\beta y \cdot x) \cdot (y\lambda)) -$
$b(\beta y \cdot x, (\beta y \cdot x) \cdot (x\mu)) - b(\beta y \cdot x, (\beta y \cdot x) \cdot (\alpha x \cdot y)) - b(\beta y \cdot x, (\beta y \cdot x) \cdot (\beta y \cdot x)) -$
$3 \mu b(x, x \cdot x) q(y) + 3 \alpha \beta b(x, y) q(x) q(y) - 4 b(x, y) q(x) q(y) - b(y, x) q(x) q(y) +$
$3 b(x, y) q(y\lambda) + 3 b(x, y) q(x\mu) + 3 b(x, y) q(\alpha x \cdot y) + 3 b(x, y) q(\beta y \cdot x)$

**Expand[%]**

$b(y, y \cdot y) \lambda^3 - 3 \mu q(y) \lambda^3 - 3 \mu^2 b(x, y) \lambda^2 + 3 \mu b(y \lambda, x \mu) \lambda + 3 \mu b(y \lambda, \alpha x \cdot y) \lambda + 3 \mu b(y \lambda, \beta y \cdot x) \lambda + 3 \mu b(x \mu, \alpha x \cdot y) \lambda + 3 \mu b(x \mu, \beta y \cdot x) \lambda - 3 \mu b(x \cdot x, y \cdot y) \lambda - 3 \mu b(x \cdot y, y \cdot y) \lambda + 3 \mu b(\alpha x \cdot y, \beta y \cdot x) \lambda - 3 \mu^3 q(x) \lambda - 3 b(y, y \cdot y) q(x) \lambda + 9 \mu q(x) q(y) \lambda + 3 \mu q(y \lambda) \lambda + 3 \mu q(x \mu) \lambda + 3 \mu q(\alpha x \cdot y) \lambda + 3 \mu q(\beta y \cdot x) \lambda + \mu^3 b(x, x \cdot x) - b(y \lambda, (y \lambda) \cdot (y \lambda)) - b(y \lambda, (y \lambda) \cdot (x \mu)) - b(y \lambda, (y \lambda) \cdot (\alpha x \cdot y)) - b(y \lambda, (y \lambda) \cdot (\beta y \cdot x)) - b(y \lambda, (x \mu) \cdot (y \lambda)) - b(y \lambda, (x \mu) \cdot (x \mu)) - b(y \lambda, (x \mu) \cdot (\alpha x \cdot y)) - b(y \lambda, (x \mu) \cdot (\beta y \cdot x)) - b(y \lambda, (\alpha x \cdot y) \cdot (y \lambda)) - b(y \lambda, (\alpha x \cdot y) \cdot (x \mu)) - b(y \lambda, (\alpha x \cdot y) \cdot (\alpha x \cdot y)) - b(y \lambda, (\alpha x \cdot y) \cdot (\beta y \cdot x)) - b(y \lambda, (\beta y \cdot x) \cdot (y \lambda)) - b(y \lambda, (\beta y \cdot x) \cdot (x \mu)) - b(y \lambda, (\beta y \cdot x) \cdot (\alpha x \cdot y)) - b(y \lambda, (\beta y \cdot x) \cdot (\beta y \cdot x)) - b(x \mu, (y \lambda) \cdot (y \lambda)) - b(x \mu, (y \lambda) \cdot (x \mu)) - b(x \mu, (y \lambda) \cdot (\alpha x \cdot y)) - b(x \mu, (y \lambda) \cdot (\beta y \cdot x)) - b(x \mu, (x \mu) \cdot (y \lambda)) - b(x \mu, (x \mu) \cdot (x \mu)) - b(x \mu, (x \mu) \cdot (\alpha x \cdot y)) - b(x \mu, (x \mu) \cdot (\beta y \cdot x)) - b(x \mu, (\alpha x \cdot y) \cdot (y \lambda)) - b(x \mu, (\alpha x \cdot y) \cdot (x \mu)) - b(x \mu, (\alpha x \cdot y) \cdot (\alpha x \cdot y)) - b(x \mu, (\alpha x \cdot y) \cdot (\beta y \cdot x)) - b(x \mu, (\beta y \cdot x) \cdot (y \lambda)) - b(x \mu, (\beta y \cdot x) \cdot (x \mu)) - b(x \mu, (\beta y \cdot x) \cdot (\alpha x \cdot y)) - b(x \mu, (\beta y \cdot x) \cdot (\beta y \cdot x)) + 3 \alpha \beta b(x \cdot y, \lambda \mu x \cdot y) - 3 b(x \cdot y, \lambda \mu x \cdot y) + b(x \cdot y, (x \cdot x) \cdot (y \cdot y)) + 3 \alpha \beta b(x \cdot y, -\lambda (x \cdot y) \cdot y) - 3 b(x \cdot y, -\lambda (x \cdot y) \cdot y) + \alpha \beta^2 b(x \cdot y, (x \cdot y) \cdot (x \cdot y)) + \alpha \beta b(x \cdot y, (x \cdot y) \cdot (x \cdot y)) - \beta b(x \cdot y, (x \cdot y) \cdot (x \cdot y)) - b(x \cdot y, (x \cdot y) \cdot (x \cdot y)) - \alpha \beta^2 b(x \cdot y, (x \cdot y) \cdot (y \cdot x)) - \alpha \beta b(x \cdot y, (x \cdot y) \cdot (y \cdot x)) + \beta b(x \cdot y, (x \cdot y) \cdot (y \cdot x)) + b(x \cdot y, (x \cdot y) \cdot (y \cdot x)) + 3 \alpha \beta b(x \cdot y, \mu (y \cdot x) \cdot x) - 3 b(x \cdot y, \mu (y \cdot x) \cdot x) - \alpha \beta^2 b(x \cdot y, (y \cdot x) \cdot (x \cdot y)) - \alpha \beta b(x \cdot y, (y \cdot x) \cdot (x \cdot y)) + \beta b(x \cdot y, (y \cdot x) \cdot (x \cdot y)) + b(x \cdot y, (y \cdot x) \cdot (x \cdot y)) + \alpha \beta^2 b(x \cdot y, (y \cdot x) \cdot (y \cdot x)) + \alpha \beta b(x \cdot y, (y \cdot x) \cdot (y \cdot x)) - \beta b(x \cdot y, (y \cdot x) \cdot (y \cdot x)) - b(x \cdot y, (y \cdot x) \cdot (y \cdot x)) - b(\alpha x \cdot y, (y \lambda) \cdot (y \lambda)) - b(\alpha x \cdot y, (y \lambda) \cdot (x \mu)) - b(\alpha x \cdot y, (y \lambda) \cdot (\alpha x \cdot y)) - b(\alpha x \cdot y, (y \lambda) \cdot (\beta y \cdot x)) - b(\alpha x \cdot y, (x \mu) \cdot (y \lambda)) - b(\alpha x \cdot y, (x \mu) \cdot (x \mu)) - b(\alpha x \cdot y, (x \mu) \cdot (\alpha x \cdot y)) - b(\alpha x \cdot y, (x \mu) \cdot (\beta y \cdot x)) - b(\alpha x \cdot y, (\alpha x \cdot y) \cdot (y \lambda)) - b(\alpha x \cdot y, (\alpha x \cdot y) \cdot (x \mu)) - b(\alpha x \cdot y, (\alpha x \cdot y) \cdot (\alpha x \cdot y)) - b(\alpha x \cdot y, (\alpha x \cdot y) \cdot (\beta y \cdot x)) - b(\alpha x \cdot y, (\beta y \cdot x) \cdot (y \lambda)) - b(\alpha x \cdot y, (\beta y \cdot x) \cdot (x \mu)) - b(\alpha x \cdot y, (\beta y \cdot x) \cdot (\alpha x \cdot y)) - b(\alpha x \cdot y, (\beta y \cdot x) \cdot (\beta y \cdot x)) + 3 b(x \cdot (y \lambda), y \cdot (x \mu)) + 3 b(x \cdot (y \lambda), y \cdot (\alpha x \cdot y)) + 3 b(x \cdot (y \lambda), y \cdot (\beta y \cdot x)) + 3 b(x \cdot (x \mu), y \cdot (y \lambda)) + 3 b(x \cdot (x \mu), y \cdot (\alpha x \cdot y)) + 3 b(x \cdot (x \mu), y \cdot (\beta y \cdot x)) - b(x \cdot (x \cdot y), y \cdot (y \cdot x)) - 3 \alpha \beta b(x \cdot (x \cdot y), (y \cdot y) \cdot x) + 3 b(x \cdot (x \cdot y), (y \cdot y) \cdot x) + 3 b(x \cdot (\alpha x \cdot y), y \cdot (y \lambda)) + 3 b(x \cdot (\alpha x \cdot y), y \cdot (x \mu)) + 3 b(x \cdot (\alpha x \cdot y), y \cdot (\beta y \cdot x)) + 3 b(x \cdot (\beta y \cdot x), y \cdot (y \lambda)) + 3 b(x \cdot (\beta y \cdot x), y \cdot (x \mu)) + 3 b(x \cdot (\beta y \cdot x), y \cdot (\alpha x \cdot y)) - b(x \cdot (y \cdot y), y \cdot (x \cdot y)) + b(x \cdot (y \cdot y), (x \cdot x) \cdot y) - 3 \alpha \beta b(y \cdot x, \lambda \mu x \cdot y) + 3 b(y \cdot x, \lambda \mu x \cdot y) - 3 \alpha \beta b(y \cdot x, -\lambda (x \cdot y) \cdot y) + 3 b(y \cdot x, -\lambda (x \cdot y) \cdot y) - \alpha \beta^2 b(y \cdot x, (x \cdot y) \cdot (x \cdot y)) - \alpha \beta b(y \cdot x, (x \cdot y) \cdot (x \cdot y)) + \beta b(y \cdot x, (x \cdot y) \cdot (x \cdot y)) + b(y \cdot x, (x \cdot y) \cdot (x \cdot y)) + \alpha \beta^2 b(y \cdot x, (x \cdot y) \cdot (y \cdot x)) + \alpha \beta b(y \cdot x, (x \cdot y) \cdot (y \cdot x)) - \beta b(y \cdot x, (x \cdot y) \cdot (y \cdot x)) - b(y \cdot x, (x \cdot y) \cdot (y \cdot x)) - 3 \alpha \beta b(y \cdot x, \mu (y \cdot x) \cdot x) + 3 b(y \cdot x, \mu (y \cdot x) \cdot x) + \alpha \beta^2 b(y \cdot x, (y \cdot x) \cdot (x \cdot y)) + \alpha \beta b(y \cdot x, (y \cdot x) \cdot (x \cdot y)) - \beta b(y \cdot x, (y \cdot x) \cdot (x \cdot y)) - b(y \cdot x, (y \cdot x) \cdot (x \cdot y)) - \alpha \beta^2 b(y \cdot x, (y \cdot x) \cdot (y \cdot x)) - \alpha \beta b(y \cdot x, (y \cdot x) \cdot (y \cdot x)) + \beta b(y \cdot x, (y \cdot x) \cdot (y \cdot x)) + b(y \cdot x, (y \cdot x) \cdot (y \cdot x)) - b(\beta y \cdot x, (y \lambda) \cdot (y \lambda)) - b(\beta y \cdot x, (y \lambda) \cdot (x \mu)) - b(\beta y \cdot x, (y \lambda) \cdot (\alpha x \cdot y)) - b(\beta y \cdot x, (y \lambda) \cdot (\beta y \cdot x)) - b(\beta y \cdot x, (x \mu) \cdot (y \lambda)) - b(\beta y \cdot x, (x \mu) \cdot (x \mu)) - b(\beta y \cdot x, (x \mu) \cdot (\alpha x \cdot y)) - b(\beta y \cdot x, (x \mu) \cdot (\beta y \cdot x)) - b(\beta y \cdot x, (\alpha x \cdot y) \cdot (y \lambda)) - b(\beta y \cdot x, (\alpha x \cdot y) \cdot (x \mu)) - b(\beta y \cdot x, (\alpha x \cdot y) \cdot (\alpha x \cdot y)) - b(\beta y \cdot x, (\alpha x \cdot y) \cdot (\beta y \cdot x)) - b(\beta y \cdot x, (\beta y \cdot x) \cdot (y \lambda)) - b(\beta y \cdot x, (\beta y \cdot x) \cdot (x \mu)) - b(\beta y \cdot x, (\beta y \cdot x) \cdot (\alpha x \cdot y)) - b(\beta y \cdot x, (\beta y \cdot x) \cdot (\beta y \cdot x)) - 3 \mu b(x, x \cdot x) q(y) + 3 \alpha \beta b(x, y) q(x) q(y) - 4 b(x, y) q(x) q(y) - b(y, x) q(x) q(y) + 3 b(x, y) q(y \lambda) + 3 b(x, y) q(x \mu) + 3 b(x, y) q(\alpha x \cdot y) + 3 b(x, y) q(\beta y \cdot x)$

**% //. assocb**

$b(y, y \cdot y) \lambda^3 - 3 \mu q(y) \lambda^3 - 3 \mu^2 b(x, y) \lambda^2 + 3 \mu b(y \lambda, x \mu) \lambda + 3 \mu b(y \lambda, \alpha x \cdot y) \lambda +$
$\quad 3 \mu b(y \lambda, \beta y \cdot x) \lambda + 3 \mu b(x \mu, \alpha x \cdot y) \lambda + 3 \mu b(x \mu, \beta y \cdot x) \lambda - 3 \mu b(x \cdot x, y \cdot y) \lambda -$
$\quad 3 \mu b(x \cdot y, y \cdot x) \lambda + 3 \mu b(\alpha x \cdot y, \beta y \cdot x) \lambda - 3 \mu^3 q(x) \lambda - 3 b(y, y \cdot y) q(x) \lambda +$
$\quad 9 \mu q(x) q(y) \lambda + 3 \mu q(y \lambda) \lambda + 3 \mu q(x \mu) \lambda + 3 \mu q(\alpha x \cdot y) \lambda + 3 \mu q(\beta y \cdot x) \lambda +$
$\quad \mu^3 b(x, x \cdot x) - b(y \lambda, (y \lambda) \cdot (y \lambda)) - b(y \lambda, (y \lambda) \cdot (x \mu)) - b(y \lambda, (y \lambda) \cdot (\alpha x \cdot y)) -$
$\quad b(y \lambda, (y \lambda) \cdot (\beta y \cdot x)) - b(y \lambda, (x \mu) \cdot (y \lambda)) - b(y \lambda, (x \mu) \cdot (x \mu)) - b(y \lambda, (x \mu) \cdot (\alpha x \cdot y)) -$
$\quad b(y \lambda, (x \mu) \cdot (\beta y \cdot x)) - b(y \lambda, (\alpha x \cdot y) \cdot (y \lambda)) - b(y \lambda, (\alpha x \cdot y) \cdot (x \mu)) -$
$\quad b(y \lambda, (\alpha x \cdot y) \cdot (\alpha x \cdot y)) - b(y \lambda, (\alpha x \cdot y) \cdot (\beta y \cdot x)) - b(y \lambda, (\beta y \cdot x) \cdot (y \lambda)) -$
$\quad b(y \lambda, (\beta y \cdot x) \cdot (x \mu)) - b(y \lambda, (\beta y \cdot x) \cdot (\alpha x \cdot y)) - b(y \lambda, (\beta y \cdot x) \cdot (\beta y \cdot x)) -$
$\quad b(x \mu, (y \lambda) \cdot (y \lambda)) - b(x \mu, (y \lambda) \cdot (x \mu)) - b(x \mu, (y \lambda) \cdot (\alpha x \cdot y)) - b(x \mu, (y \lambda) \cdot (\beta y \cdot x)) -$
$\quad b(x \mu, (x \mu) \cdot (y \lambda)) - b(x \mu, (x \mu) \cdot (x \mu)) - b(x \mu, (x \mu) \cdot (\alpha x \cdot y)) -$
$\quad b(x \mu, (x \mu) \cdot (\beta y \cdot x)) - b(x \mu, (\alpha x \cdot y) \cdot (y \lambda)) - b(x \mu, (\alpha x \cdot y) \cdot (x \mu)) -$
$\quad b(x \mu, (\alpha x \cdot y) \cdot (\alpha x \cdot y)) - b(x \mu, (\alpha x \cdot y) \cdot (\beta y \cdot x)) - b(x \mu, (\beta y \cdot x) \cdot (y \lambda)) -$
$\quad b(x \mu, (\beta y \cdot x) \cdot (x \mu)) - b(x \mu, (\beta y \cdot x) \cdot (\alpha x \cdot y)) - b(x \mu, (\beta y \cdot x) \cdot (\beta y \cdot x)) +$
$\quad 3 \alpha \beta b(x \cdot y, \lambda \mu x \cdot y) - 3 b(x \cdot y, \lambda \mu x \cdot y) + b(x \cdot y, (x \cdot x) \cdot (y \cdot y)) +$
$\quad 3 \alpha \beta b(x \cdot y, -\lambda (x \cdot y) \cdot y) - 3 b(x \cdot y, -\lambda (x \cdot y) \cdot y) + \alpha \beta^2 b(x \cdot y, (x \cdot y) \cdot (x \cdot y)) +$
$\quad \alpha \beta b(x \cdot y, (x \cdot y) \cdot (x \cdot y)) - \beta b(x \cdot y, (x \cdot y) \cdot (x \cdot y)) - b(x \cdot y, (x \cdot y) \cdot (x \cdot y)) +$
$\quad 3 \alpha \beta b(x \cdot y, \mu (y \cdot x) \cdot x) - 3 b(x \cdot y, \mu (y \cdot x) \cdot x) - b(\alpha x \cdot y, (y \lambda) \cdot (y \lambda)) -$
$\quad b(\alpha x \cdot y, (y \lambda) \cdot (x \mu)) - b(\alpha x \cdot y, (y \lambda) \cdot (\alpha x \cdot y)) - b(\alpha x \cdot y, (y \lambda) \cdot (\beta y \cdot x)) -$
$\quad b(\alpha x \cdot y, (x \mu) \cdot (y \lambda)) - b(\alpha x \cdot y, (x \mu) \cdot (x \mu)) - b(\alpha x \cdot y, (x \mu) \cdot (\alpha x \cdot y)) -$
$\quad b(\alpha x \cdot y, (x \mu) \cdot (\beta y \cdot x)) - b(\alpha x \cdot y, (\alpha x \cdot y) \cdot (y \lambda)) - b(\alpha x \cdot y, (\alpha x \cdot y) \cdot (x \mu)) -$
$\quad b(\alpha x \cdot y, (\alpha x \cdot y) \cdot (\alpha x \cdot y)) - b(\alpha x \cdot y, (\alpha x \cdot y) \cdot (\beta y \cdot x)) - b(\alpha x \cdot y, (\beta y \cdot x) \cdot (y \lambda)) -$
$\quad b(\alpha x \cdot y, (\beta y \cdot x) \cdot (x \mu)) - b(\alpha x \cdot y, (\beta y \cdot x) \cdot (\alpha x \cdot y)) - b(\alpha x \cdot y, (\beta y \cdot x) \cdot (\beta y \cdot x)) +$
$\quad 3 b(x \cdot (y \lambda), y \cdot (x \mu)) + 3 b(x \cdot (y \lambda), y \cdot (\alpha x \cdot y)) + 3 b(x \cdot (y \lambda), y \cdot (\beta y \cdot x)) +$
$\quad 3 b(x \cdot (x \mu), y \cdot (y \lambda)) + 3 b(x \cdot (x \mu), y \cdot (\alpha x \cdot y)) + 3 b(x \cdot (x \mu), y \cdot (\beta y \cdot x)) -$
$\quad b(x \cdot (x \cdot y), y \cdot (y \cdot x)) - 3 \alpha \beta b(x \cdot (x \cdot y), (y \cdot y) \cdot x) + 3 b(x \cdot (x \cdot y), (y \cdot y) \cdot x) +$
$\quad 3 b(x \cdot (\alpha x \cdot y), y \cdot (y \lambda)) + 3 b(x \cdot (\alpha x \cdot y), y \cdot (x \mu)) + 3 b(x \cdot (\alpha x \cdot y), y \cdot (\beta y \cdot x)) +$
$\quad 3 b(x \cdot (\beta y \cdot x), y \cdot (y \lambda)) + 3 b(x \cdot (\beta y \cdot x), y \cdot (x \mu)) + 3 b(x \cdot (\beta y \cdot x), y \cdot (\alpha x \cdot y)) -$
$\quad b(x \cdot (y \cdot y), y \cdot (x \cdot x)) + b(x \cdot (y \cdot y), (x \cdot x) \cdot y) - 3 \alpha \beta b(y \cdot x, \lambda \mu x \cdot y) +$
$\quad 3 b(y \cdot x, \lambda \mu x \cdot y) - 3 \alpha \beta b(y \cdot x, -\lambda (x \cdot y) \cdot y) + 3 b(y \cdot x, -\lambda (x \cdot y) \cdot y) -$
$\quad 3 \alpha \beta b(y \cdot x, \mu (y \cdot x) \cdot x) + 3 b(y \cdot x, \mu (y \cdot x) \cdot x) - \alpha \beta^2 b(y \cdot x, (y \cdot x) \cdot (y \cdot x)) -$
$\quad \alpha \beta b(y \cdot x, (y \cdot x) \cdot (y \cdot x)) + \beta b(y \cdot x, (y \cdot x) \cdot (y \cdot x)) + b(y \cdot x, (y \cdot x) \cdot (y \cdot x)) -$
$\quad b(\beta y \cdot x, (y \lambda) \cdot (y \lambda)) - b(\beta y \cdot x, (y \lambda) \cdot (x \mu)) - b(\beta y \cdot x, (y \lambda) \cdot (\alpha x \cdot y)) -$
$\quad b(\beta y \cdot x, (y \lambda) \cdot (\beta y \cdot x)) - b(\beta y \cdot x, (x \mu) \cdot (y \lambda)) - b(\beta y \cdot x, (x \mu) \cdot (x \mu)) -$
$\quad b(\beta y \cdot x, (x \mu) \cdot (\alpha x \cdot y)) - b(\beta y \cdot x, (x \mu) \cdot (\beta y \cdot x)) - b(\beta y \cdot x, (\alpha x \cdot y) \cdot (y \lambda)) -$
$\quad b(\beta y \cdot x, (\alpha x \cdot y) \cdot (x \mu)) - b(\beta y \cdot x, (\alpha x \cdot y) \cdot (\alpha x \cdot y)) - b(\beta y \cdot x, (\alpha x \cdot y) \cdot (\beta y \cdot x)) -$
$\quad b(\beta y \cdot x, (\beta y \cdot x) \cdot (y \lambda)) - b(\beta y \cdot x, (\beta y \cdot x) \cdot (x \mu)) - b(\beta y \cdot x, (\beta y \cdot x) \cdot (\alpha x \cdot y)) -$
$\quad b(\beta y \cdot x, (\beta y \cdot x) \cdot (\beta y \cdot x)) - 3 \alpha \beta^2 b((x \cdot y) \cdot (x \cdot y), y \cdot x) - 3 \alpha \beta b((x \cdot y) \cdot (x \cdot y), y \cdot x) +$
$\quad 3 \beta b((x \cdot y) \cdot (x \cdot y), y \cdot x) + 3 b((x \cdot y) \cdot (x \cdot y), y \cdot x) + 3 \alpha \beta^2 b((y \cdot x) \cdot (y \cdot x), x \cdot y) +$
$\quad 3 \alpha \beta b((y \cdot x) \cdot (y \cdot x), x \cdot y) - 3 \beta b((y \cdot x) \cdot (y \cdot x), x \cdot y) - 3 b((y \cdot x) \cdot (y \cdot x), x \cdot y) -$
$\quad 3 \mu b(x, x \cdot x) q(y) + 3 \alpha \beta b(x, y) q(x) q(y) - 4 b(x, y) q(x) q(y) - b(y, x) q(x) q(y) +$
$\quad 3 b(x, y) q(y \lambda) + 3 b(x, y) q(x \mu) + 3 b(x, y) q(\alpha x \cdot y) + 3 b(x, y) q(\beta y \cdot x)$

**% // Expand**

$b(y, y \cdot y)\lambda^3 - 3\mu q(y)\lambda^3 - 3\mu^2 b(x, y)\lambda^2 + 3\mu b(y\lambda, x\mu)\lambda + 3\mu b(y\lambda, \alpha x \cdot y)\lambda +$
$3\mu b(y\lambda, \beta y \cdot x)\lambda + 3\mu b(x\mu, \alpha x \cdot y)\lambda + 3\mu b(x\mu, \beta y \cdot x)\lambda - 3\mu b(x \cdot x, y \cdot y)\lambda -$
$3\mu b(x \cdot y, y \cdot y)\lambda + 3\mu b(\alpha x \cdot y, \beta y \cdot x)\lambda - 3\mu^3 q(x)\lambda - 3 b(y, y \cdot y) q(x)\lambda +$
$9\mu q(x) q(y)\lambda + 3\mu q(y\lambda)\lambda + 3\mu q(x\mu)\lambda + 3\mu q(\alpha x \cdot y)\lambda + 3\mu q(\beta y \cdot x)\lambda +$
$\mu^3 b(x, x \cdot x) - b(y\lambda, (y\lambda) \cdot (y\lambda)) - b(y\lambda, (y\lambda) \cdot (x\mu)) - b(y\lambda, (y\lambda) \cdot (\alpha x \cdot y)) -$
$b(y\lambda, (y\lambda) \cdot (\beta y \cdot x)) - b(y\lambda, (x\mu) \cdot (y\lambda)) - b(y\lambda, (x\mu) \cdot (x\mu)) - b(y\lambda, (x\mu) \cdot (\alpha x \cdot y)) -$
$b(y\lambda, (x\mu) \cdot (\beta y \cdot x)) - b(y\lambda, (\alpha x \cdot y) \cdot (y\lambda)) - b(y\lambda, (\alpha x \cdot y) \cdot (x\mu)) -$
$b(y\lambda, (\alpha x \cdot y) \cdot (\alpha x \cdot y)) - b(y\lambda, (\alpha x \cdot y) \cdot (\beta y \cdot x)) - b(y\lambda, (\beta y \cdot x) \cdot (y\lambda)) -$
$b(y\lambda, (\beta y \cdot x) \cdot (x\mu)) - b(y\lambda, (\beta y \cdot x) \cdot (\alpha x \cdot y)) - b(y\lambda, (\beta y \cdot x) \cdot (\beta y \cdot x)) -$
$b(x\mu, (y\lambda) \cdot (y\lambda)) - b(x\mu, (y\lambda) \cdot (x\mu)) - b(x\mu, (y\lambda) \cdot (\alpha x \cdot y)) - b(x\mu, (y\lambda) \cdot (\beta y \cdot x)) -$
$b(x\mu, (x\mu) \cdot (y\lambda)) - b(x\mu, (x\mu) \cdot (x\mu)) - b(x\mu, (x\mu) \cdot (\alpha x \cdot y)) -$
$b(x\mu, (x\mu) \cdot (\beta y \cdot x)) - b(x\mu, (\alpha x \cdot y) \cdot (y\lambda)) - b(x\mu, (\alpha x \cdot y) \cdot (x\mu)) -$
$b(x\mu, (\alpha x \cdot y) \cdot (\alpha x \cdot y)) - b(x\mu, (\alpha x \cdot y) \cdot (\beta y \cdot x)) - b(x\mu, (\beta y \cdot x) \cdot (y\lambda)) -$
$b(x\mu, (\beta y \cdot x) \cdot (x\mu)) - b(x\mu, (\beta y \cdot x) \cdot (\alpha x \cdot y)) - b(x\mu, (\beta y \cdot x) \cdot (\beta y \cdot x)) +$
$3\alpha\beta b(x \cdot y, \lambda \mu x \cdot y) - 3 b(x \cdot y, \lambda \mu x \cdot y) + b(x \cdot y, (x \cdot x) \cdot (y \cdot y)) +$
$3\alpha\beta b(x \cdot y, -\lambda(x \cdot y) \cdot y) - 3 b(x \cdot y, -\lambda(x \cdot y) \cdot y) + \alpha\beta^2 b(x \cdot y, (x \cdot y) \cdot (x \cdot y)) +$
$\alpha\beta b(x \cdot y, (x \cdot y) \cdot (x \cdot y)) - \beta b(x \cdot y, (x \cdot y) \cdot (x \cdot y)) - b(x \cdot y, (x \cdot y) \cdot (x \cdot y)) +$
$3\alpha\beta b(x \cdot y, \mu(y \cdot x) \cdot x) - 3 b(x \cdot y, \mu(y \cdot x) \cdot x) - b(\alpha x \cdot y, (y\lambda) \cdot (y\lambda)) -$
$b(\alpha x \cdot y, (y\lambda) \cdot (x\mu)) - b(\alpha x \cdot y, (y\lambda) \cdot (\alpha x \cdot y)) - b(\alpha x \cdot y, (y\lambda) \cdot (\beta y \cdot x)) -$
$b(\alpha x \cdot y, (x\mu) \cdot (y\lambda)) - b(\alpha x \cdot y, (x\mu) \cdot (x\mu)) - b(\alpha x \cdot y, (x\mu) \cdot (\alpha x \cdot y)) -$
$b(\alpha x \cdot y, (x\mu) \cdot (\beta y \cdot x)) - b(\alpha x \cdot y, (\alpha x \cdot y) \cdot (y\lambda)) - b(\alpha x \cdot y, (\alpha x \cdot y) \cdot (x\mu)) -$
$b(\alpha x \cdot y, (\alpha x \cdot y) \cdot (\alpha x \cdot y)) - b(\alpha x \cdot y, (\alpha x \cdot y) \cdot (\beta y \cdot x)) - b(\alpha x \cdot y, (\beta y \cdot x) \cdot (y\lambda)) -$
$b(\alpha x \cdot y, (\beta y \cdot x) \cdot (x\mu)) - b(\alpha x \cdot y, (\beta y \cdot x) \cdot (\alpha x \cdot y)) - b(\alpha x \cdot y, (\beta y \cdot x) \cdot (\beta y \cdot x)) +$
$3 b(x \cdot (y\lambda), y \cdot (x\mu)) + 3 b(x \cdot (y\lambda), y \cdot (\alpha x \cdot y)) + 3 b(x \cdot (y\lambda), y \cdot (\beta y \cdot x)) +$
$3 b(x \cdot (x\mu), y \cdot (y\lambda)) + 3 b(x \cdot (x\mu), y \cdot (\alpha x \cdot y)) + 3 b(x \cdot (x\mu), y \cdot (\beta y \cdot x)) -$
$b(x \cdot (x \cdot y), y \cdot (y \cdot x)) - 3\alpha\beta b(x \cdot (x \cdot y), (y \cdot y) \cdot x) + 3 b(x \cdot (x \cdot y), (y \cdot y) \cdot x) +$
$3 b(x \cdot (\alpha x \cdot y), y \cdot (y\lambda)) + 3 b(x \cdot (\alpha x \cdot y), y \cdot (x\mu)) + 3 b(x \cdot (\alpha x \cdot y), y \cdot (\beta y \cdot x)) +$
$3 b(x \cdot (\beta y \cdot x), y \cdot (y\lambda)) + 3 b(x \cdot (\beta y \cdot x), y \cdot (x\mu)) + 3 b(x \cdot (\beta y \cdot x), y \cdot (\alpha x \cdot y)) -$
$b(x \cdot (y \cdot y), y \cdot (x \cdot x)) + b(x \cdot (y \cdot y), (x \cdot x) \cdot y) - 3\alpha\beta b(y \cdot x, \lambda\mu x \cdot y) +$
$3 b(y \cdot x, \lambda\mu x \cdot y) - 3\alpha\beta b(y \cdot x, -\lambda(x \cdot y) \cdot y) + 3 b(y \cdot x, -\lambda(x \cdot y) \cdot y) -$
$3\alpha\beta b(y \cdot x, \mu(y \cdot x) \cdot x) + 3 b(y \cdot x, \mu(y \cdot x) \cdot x) - \alpha\beta^2 b(y \cdot x, (y \cdot x) \cdot (y \cdot x)) -$
$\alpha\beta b(y \cdot x, (y \cdot x) \cdot (y \cdot x)) + \beta b(y \cdot x, (y \cdot x) \cdot (y \cdot x)) + b(y \cdot x, (y \cdot x) \cdot (y \cdot x)) -$
$b(\beta y \cdot x, (y\lambda) \cdot (y\lambda)) - b(\beta y \cdot x, (y\lambda) \cdot (x\mu)) - b(\beta y \cdot x, (y\lambda) \cdot (\alpha x \cdot y)) -$
$b(\beta y \cdot x, (y\lambda) \cdot (\beta y \cdot x)) - b(\beta y \cdot x, (x\mu) \cdot (y\lambda)) - b(\beta y \cdot x, (x\mu) \cdot (x\mu)) -$
$b(\beta y \cdot x, (x\mu) \cdot (\alpha x \cdot y)) - b(\beta y \cdot x, (x\mu) \cdot (\beta y \cdot x)) - b(\beta y \cdot x, (\alpha x \cdot y) \cdot (y\lambda)) -$
$b(\beta y \cdot x, (\alpha x \cdot y) \cdot (x\mu)) - b(\beta y \cdot x, (\alpha x \cdot y) \cdot (\alpha x \cdot y)) - b(\beta y \cdot x, (\alpha x \cdot y) \cdot (\beta y \cdot x)) -$
$b(\beta y \cdot x, (\beta y \cdot x) \cdot (y\lambda)) - b(\beta y \cdot x, (\beta y \cdot x) \cdot (x\mu)) - b(\beta y \cdot x, (\beta y \cdot x) \cdot (\alpha x \cdot y)) -$
$b(\beta y \cdot x, (\beta y \cdot x) \cdot (\beta y \cdot x)) - 3\alpha\beta^2 b((x \cdot y) \cdot (x \cdot y), y \cdot x) - 3\alpha\beta b((x \cdot y) \cdot (x \cdot y), y \cdot x) +$
$3\beta b((x \cdot y) \cdot (x \cdot y), y \cdot x) + 3 b((x \cdot y) \cdot (x \cdot y), y \cdot x) + 3\alpha\beta^2 b((y \cdot x) \cdot (y \cdot x), x \cdot y) +$
$3\alpha\beta b((y \cdot x) \cdot (y \cdot x), x \cdot y) - 3\beta b((y \cdot x) \cdot (y \cdot x), x \cdot y) - 3 b((y \cdot x) \cdot (y \cdot x), x \cdot y) -$
$3\mu b(x, x \cdot x) q(y) + 3\alpha\beta b(x, y) q(x) q(y) - 4 b(x, y) q(x) q(y) - b(y, x) q(x) q(y) +$
$3 b(x, y) q(y\lambda) + 3 b(x, y) q(x\mu) + 3 b(x, y) q(\alpha x \cdot y) + 3 b(x, y) q(\beta y \cdot x)$

**SOut[%, {α, β, λ, μ}]**

$-b(x \cdot y, (x \cdot y) \cdot (x \cdot y))\, \alpha^3 - \mu\, b(x, (x \cdot y) \cdot (x \cdot y))\, \alpha^2 - \lambda\, b(y, (x \cdot y) \cdot (x \cdot y))\, \alpha^2 -$
$\quad \mu\, b(x \cdot y, x \cdot (x \cdot y))\, \alpha^2 - \lambda\, b(x \cdot y, y \cdot (x \cdot y))\, \alpha^2 - \mu\, b(x \cdot y, (x \cdot y) \cdot x)\, \alpha^2 -$
$\quad \lambda\, b(x \cdot y, (x \cdot y) \cdot y)\, \alpha^2 - \beta\, b(x \cdot y, (x \cdot y) \cdot (y \cdot x))\, \alpha^2 - \beta\, b(x \cdot y, (y \cdot x) \cdot (x \cdot y))\, \alpha^2 -$
$\quad \beta\, b(y \cdot x, (x \cdot y) \cdot (x \cdot y))\, \alpha^2 + 3\, \lambda\, \mu\, q(x \cdot y)\, \alpha^2 + 3\, b(x, y)\, q(x \cdot y)\, \alpha^2 +$
$\quad 3\, \lambda\, \mu^2\, b(x, x \cdot y)\, \alpha - \mu^2\, b(x, x \cdot (x \cdot y))\, \alpha - \lambda\, \mu\, b(x, y \cdot (x \cdot y))\, \alpha - \mu^2\, b(x, (x \cdot y) \cdot x)\, \alpha -$
$\quad \lambda\, \mu\, b(x, (x \cdot y) \cdot y)\, \alpha - \beta\, \mu\, b(x, (x \cdot y) \cdot (y \cdot x))\, \alpha - \beta\, \mu\, b(x, (y \cdot x) \cdot (x \cdot y))\, \alpha +$
$\quad 3\, \lambda^2\, \mu\, b(y, x \cdot y)\, \alpha - \lambda\, \mu\, b(y, x \cdot (x \cdot y))\, \alpha - \lambda^2\, b(y, y \cdot (x \cdot y))\, \alpha - \lambda\, \mu\, b(y, (x \cdot y) \cdot x)\, \alpha -$
$\quad \lambda^2\, b(y, (x \cdot y) \cdot y)\, \alpha - \beta\, \lambda\, b(y, (x \cdot y) \cdot (y \cdot x))\, \alpha - \beta\, \lambda\, b(y, (y \cdot x) \cdot (x \cdot y))\, \alpha +$
$\quad 3\, \mu\, b(x \cdot x, y \cdot (x \cdot y))\, \alpha - \mu^2\, b(x \cdot y, x \cdot x)\, \alpha + 3\, \beta\, \lambda\, \mu\, b(x \cdot y, x \cdot y)\, \alpha - \lambda\, \mu\, b(x \cdot y, x \cdot y)\, \alpha -$
$\quad \beta\, \mu\, b(x \cdot y, x \cdot (y \cdot x))\, \alpha + 3\, \beta\, \lambda\, \mu\, b(x \cdot y, y \cdot x)\, \alpha - \lambda\, \mu\, b(x \cdot y, y \cdot x)\, \alpha - \lambda^2\, b(x \cdot y, y \cdot y)\, \alpha +$
$\quad 3\, \lambda\, b(x \cdot y, y \cdot (x \cdot y))\, \alpha - \beta\, \lambda\, b(x \cdot y, y \cdot (y \cdot x))\, \alpha - 3\, \beta\, \lambda\, b(x \cdot y, (x \cdot y) \cdot y)\, \alpha +$
$\quad \beta^2\, b(x \cdot y, (x \cdot y) \cdot (x \cdot y))\, \alpha + \beta\, b(x \cdot y, (x \cdot y) \cdot (x \cdot y))\, \alpha + 2\, \beta\, \mu\, b(x \cdot y, (y \cdot x) \cdot x)\, \alpha -$
$\quad \beta\, \lambda\, b(x \cdot y, (y \cdot x) \cdot y)\, \alpha - \beta^2\, b(x \cdot y, (y \cdot x) \cdot (y \cdot x))\, \alpha + 3\, \mu\, b(x \cdot (x \cdot y), y \cdot x)\, \alpha +$
$\quad 3\, \lambda\, b(x \cdot (x \cdot y), y \cdot y)\, \alpha + 3\, \beta\, b(x \cdot (x \cdot y), y \cdot (y \cdot x))\, \alpha - 3\, \beta\, b(x \cdot (x \cdot y), (y \cdot y) \cdot x)\, \alpha +$
$\quad 3\, \beta\, b(x \cdot (y \cdot x), y \cdot (x \cdot y))\, \alpha - 3\, \beta\, \lambda\, \mu\, b(y \cdot x, x \cdot y)\, \alpha - \beta\, \mu\, b(y \cdot x, x \cdot (x \cdot y))\, \alpha -$
$\quad \beta\, \lambda\, b(y \cdot x, y \cdot (x \cdot y))\, \alpha - \beta\, \mu\, b(y \cdot x, (x \cdot y) \cdot x)\, \alpha + 2\, \beta\, \lambda\, b(y \cdot x, (x \cdot y) \cdot y)\, \alpha -$
$\quad \beta^2\, b(y \cdot x, (x \cdot y) \cdot (y \cdot x))\, \alpha - 3\, \beta\, \mu\, b(y \cdot x, (y \cdot x) \cdot x)\, \alpha - \beta^2\, b(y \cdot x, (y \cdot x) \cdot (x \cdot y))\, \alpha -$
$\quad \beta^2\, b(y \cdot x, (y \cdot x) \cdot (y \cdot x))\, \alpha - \beta\, b(y \cdot x, (y \cdot x) \cdot (y \cdot x))\, \alpha - 3\, \beta^2\, b((x \cdot y) \cdot (x \cdot y), y \cdot x)\, \alpha -$
$\quad 3\, \beta\, b((x \cdot y) \cdot (x \cdot y), y \cdot x)\, \alpha + 3\, \beta^2\, b((y \cdot x) \cdot (y \cdot x), x \cdot y)\, \alpha + 3\, \beta\, b((y \cdot x) \cdot (y \cdot x), x \cdot y)\, \alpha +$
$\quad 3\, \beta\, b(x, y)\, q(x)\, q(y)\, \alpha - 3\, \lambda^2\, \mu^2\, b(x, y) - \lambda\, \mu^2\, b(x, x \cdot y) - \beta\, \mu^2\, b(x, x \cdot (y \cdot x)) +$
$\quad 3\, \beta\, \lambda\, \mu^2\, b(x, y \cdot x) - \lambda\, \mu^2\, b(x, y \cdot x) - \lambda^2\, \mu\, b(x, y \cdot y) - \beta\, \lambda\, \mu\, b(x, y \cdot (y \cdot x)) -$
$\quad \beta\, \mu^2\, b(x, (y \cdot x) \cdot x) - \beta\, \lambda\, \mu\, b(x, (y \cdot x) \cdot y) - \beta^2\, \mu\, b(x, (y \cdot x) \cdot (y \cdot x)) +$
$\quad 3\, \lambda^2\, \mu^2\, b(y, x) - \lambda\, \mu^2\, b(y, x \cdot x) - \lambda^2\, \mu\, b(y, x \cdot y) - \beta\, \lambda\, \mu\, b(y, x \cdot (y \cdot x)) +$
$\quad 3\, \beta\, \lambda^2\, \mu\, b(y, y \cdot x) - \lambda^2\, \mu\, b(y, y \cdot x) - \beta\, \lambda^2\, b(y, y \cdot (y \cdot x)) - \beta\, \lambda\, \mu\, b(y, (y \cdot x) \cdot x) -$
$\quad \beta\, \lambda^2\, b(y, (y \cdot x) \cdot y) - \beta^2\, \lambda\, b(y, (y \cdot x) \cdot (y \cdot x)) + 3\, \beta\, \mu\, b(x \cdot x, y \cdot (y \cdot x)) -$
$\quad 3\, \lambda\, \mu\, b(x \cdot y, x \cdot y) + 3\, \beta\, \lambda\, b(x \cdot y, y \cdot (y \cdot x)) + b(x \cdot y, (x \cdot x) \cdot (y \cdot y)) +$
$\quad 3\, \lambda\, b(x \cdot y, (x \cdot y) \cdot y) - \beta\, b(x \cdot y, (x \cdot y) \cdot (x \cdot y)) - b(x \cdot y, (x \cdot y) \cdot (x \cdot y)) -$
$\quad 3\, \mu\, b(x \cdot y, (y \cdot x) \cdot x) - b(x \cdot (x \cdot y), y \cdot (y \cdot x)) + 3\, b(x \cdot (x \cdot y), (y \cdot y) \cdot x) +$
$\quad 3\, \beta\, \mu\, b(x \cdot (y \cdot x), y \cdot x) + 3\, \beta\, \lambda\, b(x \cdot (y \cdot x), y \cdot y) - b(x \cdot (y \cdot y), y \cdot (x \cdot x)) +$
$\quad b(x \cdot (y \cdot y), (x \cdot x) \cdot y) - \beta\, \mu^2\, b(y \cdot x, x \cdot x) - \beta\, \lambda\, \mu\, b(y \cdot x, x \cdot y) + 3\, \lambda\, \mu\, b(y \cdot x, x \cdot y) -$
$\quad \beta^2\, \mu\, b(y \cdot x, x \cdot (y \cdot x)) - \beta\, \lambda\, \mu\, b(y \cdot x, y \cdot x) - \beta\, \lambda^2\, b(y \cdot x, y \cdot y) - \beta^2\, \lambda\, b(y \cdot x, y \cdot (y \cdot x)) -$
$\quad 3\, \lambda\, b(y \cdot x, (x \cdot y) \cdot y) - \beta^2\, \mu\, b(y \cdot x, (y \cdot x) \cdot x) + 3\, \mu\, b(y \cdot x, (y \cdot x) \cdot x) -$
$\quad \beta^2\, \lambda\, b(y \cdot x, (y \cdot x) \cdot y) - \beta^3\, b(y \cdot x, (y \cdot x) \cdot (y \cdot x)) + \beta\, b(y \cdot x, (y \cdot x) \cdot (y \cdot x)) +$
$\quad b(y \cdot x, (y \cdot x) \cdot (y \cdot x)) + 3\, \beta\, b((x \cdot y) \cdot (x \cdot y), y \cdot x) + 3\, b((x \cdot y) \cdot (x \cdot y), y \cdot x) -$
$\quad 3\, \beta\, b((y \cdot x) \cdot (y \cdot x), x \cdot y) - 3\, b((y \cdot x) \cdot (y \cdot x), x \cdot y) + 3\, \mu^2\, b(x, y)\, q(x) -$
$\quad 3\, \lambda\, b(y, y \cdot y)\, q(x) + 3\, \lambda^2\, b(x, y)\, q(y) - 3\, \mu\, b(x, x \cdot x)\, q(y) + 9\, \lambda\, \mu\, q(x)\, q(y) -$
$\quad 4\, b(x, y)\, q(x)\, q(y) - b(y, x)\, q(x)\, q(y) + 3\, \beta^2\, \lambda\, \mu\, q(y \cdot x) + 3\, \beta^2\, b(x, y)\, q(y \cdot x)$

**% //. rules2**

$-b(x \cdot y, (x \cdot y) \cdot (x \cdot y)) \, \alpha^3 - \mu \, b(x, (x \cdot y) \cdot (x \cdot y)) \, \alpha^2 - \lambda \, b(y, (x \cdot y) \cdot (x \cdot y)) \, \alpha^2 - \beta \, b(x \cdot y, (x \cdot y) \cdot (y \cdot x)) \, \alpha^2 - \beta \, b(x \cdot y, (y \cdot x) \cdot (x \cdot y)) \, \alpha^2 - \beta \, b(y \cdot x, (x \cdot y) \cdot (x \cdot y)) \, \alpha^2 - \mu \, b(y, x \cdot y) \, q(x) \, \alpha^2 - \mu \, b(x \cdot y, y) \, q(x) \, \alpha^2 - \lambda \, b(x, x \cdot y) \, q(y) \, \alpha^2 - \lambda \, b(x \cdot y, x) \, q(y) \, \alpha^2 + 3 \, \lambda \, \mu \, q(x) \, q(y) \, \alpha^2 + 3 \, b(x, y) \, q(x) \, q(y) \, \alpha^2 + 3 \, \lambda \, \mu^2 \, b(x, x \cdot y) \, \alpha - \mu^2 \, b(x, x \cdot (x \cdot y)) \, \alpha - \lambda \, \mu \, b(x, (x \cdot y) \cdot y) \, \alpha - \beta \, \mu \, b(x, (x \cdot y) \cdot (y \cdot x)) \, \alpha - \beta \, \mu \, b(x, (y \cdot x) \cdot (x \cdot y)) \, \alpha + 3 \, \lambda^2 \, \mu \, b(y, x \cdot y) \, \alpha - \lambda \, \mu \, b(y, x \cdot (x \cdot y)) \, \alpha - \lambda^2 \, b(y, (x \cdot y) \cdot y) \, \alpha - \beta \, \lambda \, b(y, (x \cdot y) \cdot (y \cdot x)) \, \alpha - \beta \, \lambda \, b(y, (y \cdot x) \cdot (x \cdot y)) \, \alpha + 3 \, \beta \, \lambda \, \mu \, b(x \cdot y, y \cdot x) \, \alpha - \lambda \, \mu \, b(x \cdot y, y \cdot x) \, \alpha - \beta \, \lambda \, b(x \cdot y, y \cdot (y \cdot x)) \, \alpha + \beta^2 \, b(x \cdot y, (x \cdot y) \cdot (x \cdot y)) \, \alpha + \beta \, b(x \cdot y, (x \cdot y) \cdot (x \cdot y)) \, \alpha + 2 \, \beta \, \mu \, b(x \cdot y, (y \cdot x) \cdot x) \, \alpha - \beta^2 \, b(x \cdot y, (y \cdot x) \cdot (y \cdot x)) \, \alpha + 3 \, \mu \, b(x \cdot (x \cdot y), y \cdot x) \, \alpha + 3 \, \lambda \, b(x \cdot (x \cdot y), y \cdot y) \, \alpha + 3 \, \beta \, b(x \cdot (x \cdot y), y \cdot (y \cdot x)) \, \alpha - 3 \, \beta \, b(x \cdot (x \cdot y), (y \cdot y) \cdot x) \, \alpha - 3 \, \beta \, \lambda \, \mu \, b(y \cdot x, x \cdot y) \, \alpha - \beta \, \mu \, b(y \cdot x, x \cdot (x \cdot y)) \, \alpha + 2 \, \beta \, \lambda \, b(y \cdot x, (x \cdot y) \cdot y) \, \alpha - \beta^2 \, b(y \cdot x, (x \cdot y) \cdot (y \cdot x)) \, \alpha - \beta^2 \, b(y \cdot x, (y \cdot x) \cdot (x \cdot y)) \, \alpha - \beta^2 \, b(y \cdot x, (y \cdot x) \cdot (y \cdot x)) \, \alpha - \beta \, b(y \cdot x, (y \cdot x) \cdot (y \cdot x)) \, \alpha - 3 \, \beta^2 \, b((x \cdot y) \cdot (x \cdot y), y \cdot x) \, \alpha - 3 \, \beta \, b((x \cdot y) \cdot (x \cdot y), y \cdot x) \, \alpha + 3 \, \beta^2 \, b((y \cdot x) \cdot (y \cdot x), x \cdot y) \, \alpha + 3 \, \beta \, b((y \cdot x) \cdot (y \cdot x), x \cdot y) \, \alpha - \mu^2 \, b(x, y) \, q(x) \, \alpha - \mu^2 \, b(y, x) \, q(x) \, \alpha - \beta \, \mu \, b(y, x \cdot y) \, q(x) \, \alpha - 4 \, \beta \, \mu \, b(y, y \cdot x) \, q(x) \, \alpha - \lambda^2 \, b(x, y) \, q(y) \, \alpha - 4 \, \beta \, \lambda \, b(x, x \cdot y) \, q(y) \, \alpha - \beta \, \lambda \, b(x, y \cdot x) \, q(y) \, \alpha - \lambda^2 \, b(y, x) \, q(y) \, \alpha + 3 \, \mu \, b(x \cdot x, x) \, q(y) \, \alpha + 3 \, \lambda \, b(x \cdot y, x) \, q(y) \, \alpha + 6 \, \beta \, \lambda \, \mu \, q(x) \, q(y) \, \alpha - 6 \, \lambda \, \mu \, q(x) \, q(y) \, \alpha + 3 \, \beta \, b(x, y) \, q(x) \, q(y) \, \alpha + 3 \, \beta \, b(y, x) \, q(x) \, q(y) \, \alpha - 3 \, \lambda^2 \, \mu^2 \, b(x, y) - \lambda \, \mu^2 \, b(x, x \cdot y) + 3 \, \beta \, \lambda \, \mu^2 \, b(x, y \cdot x) - \lambda \, \mu^2 \, b(x, y \cdot x) - \lambda^2 \, \mu \, b(x, y \cdot y) - \beta \, \lambda \, \mu \, b(x, y \cdot (y \cdot x)) - \beta \, \mu^2 \, b(x, (y \cdot x) \cdot x) - \beta^2 \, \mu \, b(x, (y \cdot x) \cdot (y \cdot x)) + 3 \, \lambda^2 \, \mu^2 \, b(y, x) - \lambda \, \mu^2 \, b(y, x \cdot x) - \lambda^2 \, \mu \, b(y, x \cdot y) + 3 \, \beta \, \lambda^2 \, \mu \, b(y, y \cdot x) - \lambda^2 \, \mu \, b(y, y \cdot x) - \beta \, \lambda^2 \, b(y, y \cdot (y \cdot x)) - \beta \, \lambda \, \mu \, b(y, (y \cdot x) \cdot x) - \beta^2 \, \lambda \, b(y, (y \cdot x) \cdot (y \cdot x)) + 3 \, \beta \, \mu \, b(x \cdot x, y \cdot (y \cdot x)) + 3 \, \beta \, \lambda \, b(x \cdot y, y \cdot (y \cdot x)) + b(x \cdot y, (x \cdot x) \cdot (y \cdot y)) - \beta \, b(x \cdot y, (x \cdot y) \cdot (x \cdot y)) - b(x \cdot y, (x \cdot y) \cdot (x \cdot y)) - 3 \, \mu \, b(x \cdot y, (y \cdot x) \cdot x) - b(x \cdot (x \cdot y), y \cdot (y \cdot x)) + 3 \, b(x \cdot (x \cdot y), (y \cdot y) \cdot x) - b(x \cdot (y \cdot y), y \cdot (x \cdot x)) + b(x \cdot (y \cdot y), (x \cdot x) \cdot y) - \beta \, \lambda \, \mu \, b(y \cdot x, x \cdot y) + 3 \, \lambda \, \mu \, b(y \cdot x, x \cdot y) - 3 \, \lambda \, b(y \cdot x, (x \cdot y) \cdot y) - \beta^3 \, b(y \cdot x, (y \cdot x) \cdot (y \cdot x)) + \beta \, b(y \cdot x, (y \cdot x) \cdot (y \cdot x)) + b(y \cdot x, (y \cdot x) \cdot (y \cdot x)) + 3 \, \beta \, b((x \cdot y) \cdot (x \cdot y), y \cdot x) + 3 \, b((x \cdot y) \cdot (x \cdot y), y \cdot x) - 3 \, \beta \, b((y \cdot x) \cdot (y \cdot x), x \cdot y) - 3 \, b((y \cdot x) \cdot (y \cdot x), x \cdot y) - \beta \, \mu^2 \, b(x, y) \, q(x) + 3 \, \mu^2 \, b(x, y) \, q(x) - \beta \, \mu^2 \, b(y, x) \, q(x) - \beta^2 \, \mu \, b(y, y \cdot x) \, q(x) + 3 \, \beta \, \mu \, b(y, y \cdot x) \, q(x) + 3 \, \mu \, b(y, y \cdot x) \, q(x) + 3 \, \beta \, \lambda \, b(y, y \cdot y) \, q(x) - 3 \, \lambda \, b(y, y \cdot y) \, q(x) - \beta^2 \, \mu \, b(y \cdot x, y) \, q(x) - \beta \, \lambda^2 \, b(x, y) \, q(y) + 3 \, \lambda^2 \, b(x, y) \, q(y) - 3 \, \mu \, b(x, x \cdot x) \, q(y) + 3 \, \lambda \, b(x, x \cdot y) \, q(y) - \beta^2 \, \lambda \, b(x, y \cdot x) \, q(y) - \beta \, \lambda^2 \, b(y, x) \, q(y) - \beta^2 \, \lambda \, b(y \cdot x, x) \, q(y) + 3 \, \beta^2 \, \lambda \, \mu \, q(x) \, q(y) - 6 \, \beta \, \lambda \, \mu \, q(x) \, q(y) + 3 \, \lambda \, \mu \, q(x) \, q(y) + 3 \, \beta^2 \, b(x, y) \, q(x) \, q(y) - 4 \, b(x, y) \, q(x) \, q(y) - b(y, x) \, q(x) \, q(y)$

```
% //. b[y, x] → b[x, y]
```

$-b(x \cdot y, (x \cdot y) \cdot (x \cdot y)) \alpha^3 - \mu b(x, (x \cdot y) \cdot (x \cdot y)) \alpha^2 - \lambda b(y, (x \cdot y) \cdot (x \cdot y)) \alpha^2 - \beta b(x \cdot y, (x \cdot y) \cdot (y \cdot x)) \alpha^2 - \beta b(x \cdot y, (y \cdot x) \cdot (x \cdot y)) \alpha^2 - \beta b(y \cdot x, (x \cdot y) \cdot (x \cdot y)) \alpha^2 - \mu b(y, x \cdot y) q(x) \alpha^2 - \mu b(x \cdot y, y) q(x) \alpha^2 - \lambda b(x, x \cdot y) q(y) \alpha^2 - \lambda b(x \cdot y, x) q(y) \alpha^2 + 3 \lambda \mu q(x) q(y) \alpha^2 + 3 b(x, y) q(x) q(y) \alpha^2 + 3 \lambda \mu^2 b(x, x \cdot y) \alpha - \mu^2 b(x, x \cdot (x \cdot y)) \alpha - \lambda \mu b(x, (x \cdot y) \cdot y) \alpha - \beta \mu b(x, (x \cdot y) \cdot (y \cdot x)) \alpha - \beta \mu b(x, (y \cdot x) \cdot (x \cdot y)) \alpha + 3 \lambda^2 \mu b(y, x \cdot y) \alpha - \lambda \mu b(y, x \cdot (x \cdot y)) \alpha - \lambda^2 b(y, (x \cdot y) \cdot y) \alpha - \beta \lambda b(y, (x \cdot y) \cdot (y \cdot x)) \alpha - \beta \lambda b(y, (y \cdot x) \cdot (x \cdot y)) \alpha + 3 \beta \lambda \mu b(x \cdot y, y \cdot x) \alpha - \lambda \mu b(x \cdot y, y \cdot x) \alpha - \beta \lambda b(x \cdot y, y \cdot (y \cdot x)) \alpha + \beta^2 b(x \cdot y, (x \cdot y) \cdot (x \cdot y)) + \beta b(x \cdot y, (x \cdot y) \cdot (x \cdot y)) \alpha + 2 \beta \mu b(x \cdot y, (y \cdot x) \cdot x) \alpha - \beta^2 b(x \cdot y, (y \cdot x) \cdot (y \cdot x)) \alpha + 3 \mu b(x \cdot (x \cdot y), y \cdot x) \alpha + 3 \lambda b(x \cdot (x \cdot y), y \cdot y) \alpha + 3 \beta b(x \cdot (x \cdot y), y \cdot (y \cdot x)) \alpha - 3 \beta b(x \cdot (x \cdot y), (y \cdot y) \cdot x) \alpha - 3 \beta \lambda \mu b(y \cdot x, x \cdot y) \alpha - \beta \mu b(y \cdot x, x \cdot (x \cdot y)) \alpha + 2 \beta \lambda b(y \cdot x, (x \cdot y) \cdot y) \alpha - \beta^2 b(y \cdot x, (x \cdot y) \cdot (y \cdot x)) \alpha - \beta^2 b(y \cdot x, (y \cdot x) \cdot (x \cdot y)) \alpha - \beta^2 b(y \cdot x, (y \cdot x) \cdot (y \cdot x)) \alpha - \beta b(y \cdot x, (y \cdot x) \cdot (y \cdot x)) \alpha - 3 \beta^2 b((x \cdot y) \cdot (x \cdot y), y \cdot x) \alpha - 3 \beta b((x \cdot y) \cdot (x \cdot y), y \cdot x) \alpha + 3 \beta^2 b((y \cdot x) \cdot (y \cdot x), x \cdot y) \alpha + 3 \beta b((y \cdot x) \cdot (y \cdot x), x \cdot y) \alpha - 2 \mu^2 b(x, y) q(x) \alpha - \beta \mu b(y, x \cdot y) q(x) \alpha - 4 \beta \mu b(y, y \cdot x) q(x) \alpha - 2 \lambda^2 b(x, y) q(y) \alpha - 4 \beta \lambda b(x, x \cdot y) q(y) \alpha - \beta \lambda b(x, y \cdot x) q(y) \alpha + 3 \mu b(x \cdot x, x) q(y) \alpha + 3 \lambda b(x \cdot y, x) q(y) \alpha + 6 \beta \lambda \mu q(x) q(y) \alpha - 6 \lambda \mu q(x) q(y) \alpha + 6 \beta b(x, y) q(x) q(y) \alpha - \lambda \mu^2 b(x, x \cdot y) + 3 \beta \lambda \mu^2 b(x, y \cdot x) - \lambda \mu^2 b(x, y \cdot x) - \lambda^2 \mu b(x, y \cdot y) - \beta \lambda \mu b(x, y \cdot (y \cdot x)) - \beta \mu^2 b(x, (y \cdot x) \cdot x) - \beta^2 \mu b(x, (y \cdot x) \cdot (y \cdot x)) - \lambda \mu^2 b(y, x \cdot x) - \lambda^2 \mu b(y, x \cdot y) + 3 \beta \lambda^2 \mu b(y, y \cdot x) - \lambda^2 \mu b(y, y \cdot x) - \beta \lambda^2 b(y, y \cdot (y \cdot x)) - \beta \lambda \mu b(y, (y \cdot x) \cdot x) - \beta^2 \lambda b(y, (y \cdot x) \cdot (y \cdot x)) + 3 \beta \mu b(x \cdot x, y \cdot (y \cdot x)) + 3 \beta \lambda b(x \cdot y, y \cdot (y \cdot x)) + b(x \cdot y, (x \cdot x) \cdot (y \cdot y)) - \beta b(x \cdot y, (x \cdot y) \cdot (x \cdot y)) - b(x \cdot y, (x \cdot y) \cdot (x \cdot y)) - 3 \mu b(x \cdot y, (y \cdot x) \cdot x) - b(x \cdot (x \cdot y), y \cdot (y \cdot x)) + 3 b(x \cdot (x \cdot y), (y \cdot y) \cdot x) - b(x \cdot (y \cdot y), y \cdot (x \cdot x)) + b(x \cdot (y \cdot y), (x \cdot x) \cdot y) - \beta \lambda \mu b(y \cdot x, x \cdot y) + 3 \lambda \mu b(y \cdot x, x \cdot y) - 3 \lambda b(y \cdot x, (x \cdot y) \cdot y) - \beta^3 b(y \cdot x, (y \cdot x) \cdot (y \cdot x)) + \beta b(y \cdot x, (y \cdot x) \cdot (y \cdot x)) + b(y \cdot x, (y \cdot x) \cdot (y \cdot x)) + 3 \beta b((x \cdot y) \cdot (x \cdot y), y \cdot x) + 3 b((x \cdot y) \cdot (x \cdot y), y \cdot x) - 3 \beta b((y \cdot x) \cdot (y \cdot x), x \cdot y) - 3 b((y \cdot x) \cdot (y \cdot x), x \cdot y) - 2 \beta \mu^2 b(x, y) q(x) + 3 \mu^2 b(x, y) q(x) - \beta^2 \mu b(y, y \cdot x) q(x) + 3 \beta \mu b(y, y \cdot x) q(x) + 3 \mu b(y, y \cdot x) q(x) + 3 \beta \lambda b(y, y \cdot y) q(x) - 3 \lambda b(y, y \cdot y) q(x) - \beta^2 \mu b(y \cdot x, y) q(x) - 2 \beta \lambda^2 b(x, y) q(y) + 3 \lambda^2 b(x, y) q(y) - 3 \mu b(x, x \cdot x) q(y) + 3 \lambda b(x, x \cdot y) q(y) - \beta^2 \lambda b(x, y \cdot x) q(y) - \beta^2 \lambda b(y \cdot x, x) q(y) + 3 \beta^2 \lambda \mu q(x) q(y) - 6 \beta \lambda \mu q(x) q(y) + 3 \lambda \mu q(x) q(y) + 3 \beta^2 b(x, y) q(x) q(y) - 5 b(x, y) q(x) q(y)$

```
res = %;

move1 = {b[y_, x_ · y_] → b[x, y · y], b[y_, y_ · x_] → b[x, y · y]};

Factor[Coefficient[res, λ μ²] //. move1]
```

$3 (\alpha + \beta - 1) b(y, x \cdot x)$

```
Factor[Coefficient[res, λ² μ] //. move1]
```

$3 (\alpha + \beta - 1) b(x, y \cdot y)$

```
move2 = {
   b[x, (x · y) · y] → b[x · y, y · x], b[y · x, x · y] → b[x · y, y · x],
   b[x, y · (y · x)] → b[x · y, y · x],
   b[y, x · (x · y)] → b[x · y, y · x],
   b[y, (y · x) · x] → b[x · y, y · x]};
```

```
Factor[Coefficient[res, λ μ] //. move2]
```

$-3 (\alpha + \beta - 1) (b(x \cdot y, y \cdot x) - \alpha \, q(x) \, q(y) - \beta \, q(x) \, q(y) + q(x) \, q(y))$

```
move3 = {b[y_, y_ · (y_ · x_)] → b[y · x, y · y],
   b[y_, (x_ · y_) · y_] → b[x · y, y · y]};
```

```
Coefficient[res, λ²] /. μ → 0
```

$-\beta \, b(y, y \cdot (y \cdot x)) - \alpha \, b(y, (x \cdot y) \cdot y) - 2 \alpha \, b(x, y) \, q(y) - 2 \beta \, b(x, y) \, q(y) + 3 \, b(x, y) \, q(y)$

```
% //. move3
```

$-\alpha \, b(x \cdot y, y \cdot y) - \beta \, b(y \cdot x, y \cdot y) - 2 \alpha \, b(x, y) \, q(y) - 2 \beta \, b(x, y) \, q(y) + 3 \, b(x, y) \, q(y)$

```
Factor[% //. rules1]
```

$-3 (\alpha + \beta - 1) \, b(x, y) \, q(y)$

```
Coefficient[res, μ²] /. λ → 0
```

$-\alpha \, b(x, x \cdot (x \cdot y)) - \beta \, b(x, (y \cdot x) \cdot x) - 2 \alpha \, b(x, y) \, q(x) - 2 \beta \, b(x, y) \, q(x) + 3 \, b(x, y) \, q(x)$

```
% //. move3
```

$-\alpha \, b(x \cdot y, x \cdot x) - \beta \, b(y \cdot x, x \cdot x) - 2 \alpha \, b(x, y) \, q(x) - 2 \beta \, b(x, y) \, q(x) + 3 \, b(x, y) \, q(x)$

```
% //. rules1
```

$-2 \alpha \, b(x, y) \, q(x) - 2 \beta \, b(x, y) \, q(x) + 3 \, b(x, y) \, q(x) - \alpha \, b(y, x) \, q(x) - \beta \, b(y, x) \, q(x)$

```
% /. {b[y, x] → b[x, y]}
```

$-3 \alpha \, b(x, y) \, q(x) - 3 \beta \, b(x, y) \, q(x) + 3 \, b(x, y) \, q(x)$

```
Factor[%]
```

$-3 (\alpha + \beta - 1) \, b(x, y) \, q(x)$

```
move4 = {b[x · y, x] → b[x, x · y],
   b[x, y · x] → b[x, x · y], b[y · x, x] → b[x, x · y],
   b[x · y, y · (y · x)] → b[y, (y · x) · (x · y)],
   b[y · x, (x · y) · y] → b[y, (y · x) · (x · y)],
   b[x · y, y · (y · x)] → b[y, (y · x) · (x · y)],
   b[y · x, (x · y) · y] → b[y, (y · x) · (x · y)],
   b[y, (x · y) · (x · y)] → q[y] b[x, x · y]
    (*flex sx=xt=q(x) y, ty=ys=q(y) x*),
   b[y, (y · x) · (y · x)] → q[y] b[y · x, x](*flex*),
   b[y, (x · y) · (y · x)] → q[y] b[x, y · x](*flex*)};
```

```
Coefficient[res, λ] /. μ → 0
```

$-b(y, (x \cdot y) \cdot (x \cdot y)) \alpha^2 - b(x, x \cdot y) q(y) \alpha^2 - b(x \cdot y, x) q(y) \alpha^2 - \beta b(y, (x \cdot y) \cdot (y \cdot x)) \alpha - \beta b(y, (y \cdot x) \cdot (x \cdot y)) \alpha - \beta b(x \cdot y, y \cdot (y \cdot x)) \alpha + 3 b(x \cdot (x \cdot y), y \cdot y) \alpha + 2 \beta b(y \cdot x, (x \cdot y) \cdot y) \alpha - 4 \beta b(x, x \cdot y) q(y) \alpha - \beta b(x, y \cdot x) q(y) \alpha + 3 b(x \cdot y, x) q(y) \alpha - \beta^2 b(y, (y \cdot x) \cdot (y \cdot x)) + 3 \beta b(x \cdot y, y \cdot (y \cdot x)) - 3 b(y \cdot x, (x \cdot y) \cdot y) + 3 \beta b(y, y \cdot y) q(x) - 3 b(y, y \cdot y) q(x) + 3 b(x, x \cdot y) q(y) - \beta^2 b(x, y \cdot x) q(y) - \beta^2 b(y \cdot x, x) q(y)$

```
% //. move4
```

$-3 b(x, x \cdot y) q(y) \alpha^2 + 3 b(x \cdot (x \cdot y), y \cdot y) \alpha - 6 \beta b(x, x \cdot y) q(y) \alpha + 3 b(x, x \cdot y) q(y) \alpha + 3 \beta b(y, (y \cdot x) \cdot (x \cdot y)) - 3 b(y, (y \cdot x) \cdot (x \cdot y)) + 3 \beta b(y, y \cdot y) q(x) - 3 b(y, y \cdot y) q(x) - 3 \beta^2 b(x, x \cdot y) q(y) + 3 b(x, x \cdot y) q(y)$

```
% //. β → 1 - α
```

$-3 b(x, x \cdot y) q(y) (1-\alpha)^2 + 3 b(y, (y \cdot x) \cdot (x \cdot y)) (1-\alpha) + 3 b(y, y \cdot y) q(x) (1-\alpha) - 6 \alpha b(x, x \cdot y) q(y) (1-\alpha) - 3 b(y, (y \cdot x) \cdot (x \cdot y)) + 3 \alpha b(x \cdot (x \cdot y), y \cdot y) - 3 b(y, y \cdot y) q(x) - 3 \alpha^2 b(x, x \cdot y) q(y) + 3 \alpha b(x, x \cdot y) q(y) + 3 b(x, x \cdot y) q(y)$

```
Expand[%]
```

$-3 \alpha b(y, (y \cdot x) \cdot (x \cdot y)) + 3 \alpha b(x \cdot (x \cdot y), y \cdot y) - 3 \alpha b(y, y \cdot y) q(x) + 3 \alpha b(x, x \cdot y) q(y)$

paper *

```
move5 = {b[y, x · y] → b[y, y · x],
   b[x · y, y] → b[y, y · x], b[y · x, y] → b[y, y · x],
   b[x · y, (y · x) · x] → b[x, (x · y) · (y · x)],
   b[x · (x · y), y · x] → b[x, (x · y) · (y · x)],
   b[y · x, x · (x · y)] → b[x, (x · y) · (y · x)],
   b[x, (y · x) · (y · x)] → q[x] b[y, y · x](*flex*),
   b[x, (x · y) · (x · y)] → q[x] b[x · y, y](*flex*),
   b[x, (y · x) · (x · y)] → q[x] b[y, x · y](*flex*)};
```

```
Coefficient[res, μ] /. λ → 0
```

$-b(x, (x \cdot y) \cdot (x \cdot y)) \alpha^2 - b(y, x \cdot y) q(x) \alpha^2 - b(x \cdot y, y) q(x) \alpha^2 - \beta b(x, (x \cdot y) \cdot (y \cdot x)) \alpha - \beta b(x, (y \cdot x) \cdot (x \cdot y)) \alpha + 2 \beta b(x \cdot y, (y \cdot x) \cdot x) \alpha + 3 b(x \cdot (x \cdot y), y \cdot x) \alpha - \beta b(y \cdot x, x \cdot (x \cdot y)) \alpha - \beta b(x, x \cdot y) q(x) \alpha - 4 \beta b(y, y \cdot x) q(x) \alpha + 3 b(x \cdot x, x) q(y) \alpha - \beta^2 b(x, (y \cdot x) \cdot (y \cdot x)) + 3 \beta b(x \cdot x, y \cdot (y \cdot x)) - 3 b(x \cdot y, (y \cdot x) \cdot x) - \beta^2 b(y, y \cdot x) q(x) + 3 \beta b(y, y \cdot x) q(x) + 3 b(y, y \cdot x) q(x) - \beta^2 b(y \cdot x, y) q(x) - 3 b(x, x \cdot x) q(y)$

```
% //. move5
```

$-3 b(y, y \cdot x) q(x) \alpha^2 + 3 b(x, (x \cdot y) \cdot (y \cdot x)) \alpha - 6 \beta b(y, y \cdot x) q(x) \alpha + 3 b(x \cdot x, x) q(y) \alpha - 3 b(x, (x \cdot y) \cdot (y \cdot x)) + 3 \beta b(x \cdot x, y \cdot (y \cdot x)) - 3 \beta^2 b(y, y \cdot x) q(x) + 3 \beta b(y, y \cdot x) q(x) + 3 b(y, y \cdot x) q(x) - 3 b(x, x \cdot x) q(y)$

```
% //. β → 1 - α
```

$-3\,b(y, y \cdot x)\,q(x)\,(1-\alpha)^2 + 3\,b(x \cdot x, y \cdot (y \cdot x))\,(1-\alpha) - 6\,\alpha\,b(y, y \cdot x)\,q(x)\,(1-\alpha) + 3\,b(y, y \cdot x)\,q(x)\,(1-\alpha) + 3\,\alpha\,b(x, (x \cdot y) \cdot (y \cdot x)) - 3\,b(x, (x \cdot y) \cdot (y \cdot x)) - 3\,\alpha^2\,b(y, y \cdot x)\,q(x) + 3\,b(y, y \cdot x)\,q(x) - 3\,b(x, x \cdot x)\,q(y) + 3\,\alpha\,b(x \cdot x, x)\,q(y)$

```
Expand[%]
```

$3\,\alpha\,b(x, (x \cdot y) \cdot (y \cdot x)) - 3\,b(x, (x \cdot y) \cdot (y \cdot x)) - 3\,\alpha\,b(x \cdot x, y \cdot (y \cdot x)) + 3\,b(x \cdot x, y \cdot (y \cdot x)) - 3\,\alpha\,b(y, y \cdot x)\,q(x) + 3\,b(y, y \cdot x)\,q(x) - 3\,b(x, x \cdot x)\,q(y) + 3\,\alpha\,b(x \cdot x, x)\,q(y)$

paper **

```
res /. {λ → 0, μ → 0}
```

$-b(x \cdot y, (x \cdot y) \cdot (x \cdot y))\,\alpha^3 - \beta\,b(x \cdot y, (x \cdot y) \cdot (y \cdot x))\,\alpha^2 - \beta\,b(x \cdot y, (y \cdot x) \cdot (x \cdot y))\,\alpha^2 - \beta\,b(y \cdot x, (x \cdot y) \cdot (x \cdot y))\,\alpha^2 + 3\,b(x, y)\,q(x)\,q(y)\,\alpha^2 + \beta^2\,b(x \cdot y, (x \cdot y) \cdot (x \cdot y))\,\alpha + \beta\,b(x \cdot y, (x \cdot y) \cdot (x \cdot y))\,\alpha - \beta^2\,b(x \cdot y, (y \cdot x) \cdot (y \cdot x))\,\alpha + 3\,\beta\,b(x \cdot (x \cdot y), y \cdot (y \cdot x))\,\alpha - 3\,\beta\,b(x \cdot (x \cdot y), (y \cdot y) \cdot x)\,\alpha - \beta^2\,b(y \cdot x, (x \cdot y) \cdot (y \cdot x))\,\alpha - \beta^2\,b(y \cdot x, (y \cdot x) \cdot (x \cdot y))\,\alpha - \beta^2\,b(y \cdot x, (y \cdot x) \cdot (y \cdot x))\,\alpha - \beta\,b(y \cdot x, (y \cdot x) \cdot (y \cdot x))\,\alpha - 3\,\beta^2\,b((x \cdot y) \cdot (x \cdot y), y \cdot x)\,\alpha - 3\,\beta\,b((x \cdot y) \cdot (x \cdot y), y \cdot x)\,\alpha + 3\,\beta^2\,b((y \cdot x) \cdot (y \cdot x), x \cdot y)\,\alpha + 3\,\beta\,b((y \cdot x) \cdot (y \cdot x), x \cdot y)\,\alpha + 6\,\beta\,b(x, y)\,q(x)\,q(y)\,\alpha + b(x \cdot y, (x \cdot x) \cdot (y \cdot y)) - \beta\,b(x \cdot y, (x \cdot y) \cdot (x \cdot y)) - b(x \cdot y, (x \cdot y) \cdot (x \cdot y)) - b(x \cdot (x \cdot y), y \cdot (y \cdot x)) + 3\,b(x \cdot (x \cdot y), (y \cdot y) \cdot x) - b(x \cdot (y \cdot y), y \cdot (x \cdot x)) + b(x \cdot (y \cdot y), (x \cdot x) \cdot y) - \beta^3\,b(y \cdot x, (y \cdot x) \cdot (y \cdot x)) + \beta\,b(y \cdot x, (y \cdot x) \cdot (y \cdot x)) + b(y \cdot x, (y \cdot x) \cdot (y \cdot x)) + 3\,\beta\,b((x \cdot y) \cdot (x \cdot y), y \cdot x) + 3\,b((x \cdot y) \cdot (x \cdot y), y \cdot x) - 3\,\beta\,b((y \cdot x) \cdot (y \cdot x), x \cdot y) - 3\,b((y \cdot x) \cdot (y \cdot x), x \cdot y) + 3\,\beta^2\,b(x, y)\,q(x)\,q(y) - 5\,b(x, y)\,q(x)\,q(y)$

```
Expand[% //. β → 1 - α]
```

$b(x \cdot y, (x \cdot y) \cdot (y \cdot x))\,\alpha^3 + b(x \cdot y, (y \cdot x) \cdot (x \cdot y))\,\alpha^3 - b(x \cdot y, (y \cdot x) \cdot (y \cdot x))\,\alpha^3 + b(y \cdot x, (x \cdot y) \cdot (x \cdot y))\,\alpha^3 - b(y \cdot x, (x \cdot y) \cdot (y \cdot x))\,\alpha^3 - b(y \cdot x, (y \cdot x) \cdot (x \cdot y))\,\alpha^3 - 3\,b((x \cdot y) \cdot (x \cdot y), y \cdot x)\,\alpha^3 + 3\,b((y \cdot x) \cdot (y \cdot x), x \cdot y)\,\alpha^3 - 3\,b(x \cdot y, (x \cdot y) \cdot (x \cdot y))\,\alpha^2 - b(x \cdot y, (x \cdot y) \cdot (y \cdot x))\,\alpha^2 - b(x \cdot y, (y \cdot x) \cdot (x \cdot y))\,\alpha^2 + 2\,b(x \cdot y, (y \cdot x) \cdot (y \cdot x))\,\alpha^2 - 3\,b(x \cdot (x \cdot y), y \cdot (y \cdot x))\,\alpha^2 + 3\,b(x \cdot (x \cdot y), (y \cdot y) \cdot x)\,\alpha^2 - b(y \cdot x, (x \cdot y) \cdot (x \cdot y))\,\alpha^2 + 2\,b(y \cdot x, (x \cdot y) \cdot (y \cdot x))\,\alpha^2 + 2\,b(y \cdot x, (y \cdot x) \cdot (x \cdot y))\,\alpha^2 + 9\,b((x \cdot y) \cdot (x \cdot y), y \cdot x)\,\alpha^2 - 9\,b((y \cdot x) \cdot (y \cdot x), x \cdot y)\,\alpha^2 + 3\,b(x \cdot y, (x \cdot y) \cdot (x \cdot y))\,\alpha - b(x \cdot y, (y \cdot x) \cdot (y \cdot x))\,\alpha + 3\,b(x \cdot (x \cdot y), y \cdot (y \cdot x))\,\alpha - 3\,b(x \cdot (x \cdot y), (y \cdot y) \cdot x)\,\alpha - b(y \cdot x, (x \cdot y) \cdot (y \cdot x))\,\alpha - b(y \cdot x, (y \cdot x) \cdot (x \cdot y))\,\alpha - 9\,b((x \cdot y) \cdot (x \cdot y), y \cdot x)\,\alpha + 9\,b((y \cdot x) \cdot (y \cdot x), x \cdot y)\,\alpha + b(x \cdot y, (x \cdot x) \cdot (y \cdot y)) - 2\,b(x \cdot y, (x \cdot y) \cdot (x \cdot y)) - b(x \cdot (x \cdot y), y \cdot (y \cdot x)) + 3\,b(x \cdot (x \cdot y), (y \cdot y) \cdot x) - b(x \cdot (y \cdot y), y \cdot (x \cdot x)) + b(x \cdot (y \cdot y), (x \cdot x) \cdot y) + b(y \cdot x, (y \cdot x) \cdot (y \cdot x)) + 6\,b((x \cdot y) \cdot (x \cdot y), y \cdot x) - 6\,b((y \cdot x) \cdot (y \cdot x), x \cdot y) - 2\,b(x, y)\,q(x)\,q(y)$

```
pol = %;
```

```
one = Coefficient[pol, α^2] //. assocb
```

$-3\,b(x \cdot y, (x \cdot y) \cdot (x \cdot y)) - 3\,b(x \cdot (x \cdot y), y \cdot (y \cdot x)) + 3\,b(x \cdot (x \cdot y), (y \cdot y) \cdot x) + 6\,b((x \cdot y) \cdot (x \cdot y), y \cdot x) - 3\,b((y \cdot x) \cdot (y \cdot x), x \cdot y)$

paper ***

```
two = Coefficient[pol, α] //. assocb
```

$3\,b(x \cdot y, (x \cdot y) \cdot (x \cdot y)) + 3\,b(x \cdot (x \cdot y), y \cdot (y \cdot x)) -$
$3\,b(x \cdot (x \cdot y), (y \cdot y) \cdot x) - 9\,b((x \cdot y) \cdot (x \cdot y), y \cdot x) + 6\,b((y \cdot x) \cdot (y \cdot x), x \cdot y)$

paper***

```
three = pol /. α → 0 //. assocb
```

$b(x \cdot y, (x \cdot x) \cdot (y \cdot y)) - 2\,b(x \cdot y, (x \cdot y) \cdot (x \cdot y)) -$
$b(x \cdot (x \cdot y), y \cdot (y \cdot x)) + 3\,b(x \cdot (x \cdot y), (y \cdot y) \cdot x) -$
$b(x \cdot (y \cdot y), y \cdot (x \cdot x)) + b(x \cdot (y \cdot y), (x \cdot x) \cdot y) + b(y \cdot x, (y \cdot x) \cdot (y \cdot x)) +$
$6\,b((x \cdot y) \cdot (x \cdot y), y \cdot x) - 6\,b((y \cdot x) \cdot (y \cdot x), x \cdot y) - 2\,b(x, y)\,q(x)\,q(y)$

paper ****

# ◼ References